\magnification=\magstep1
\hsize=16truecm
 
\input amstex
\TagsOnRight
\parindent=20pt
\parskip=2.8pt plus 1.3pt
\define\({\left(}
\define\){\right)}
\define\[{\left[}
\define\]{\right]}
\define\e{\varepsilon}
\define\oo{\omega}
\define\const{\text{\rm const.}\,}
\define\Sym{\text{\rm Sym}\,}
\define\supp {\sup\limits}
\define\inff{\inf\limits}
\define\summ{\sum\limits}
\define\prodd{\prod\limits}
\define\limm{\lim\limits}

\define\bigcupp{\bigcup\limits}
 
\define\sumn{\operatornamewithlimits{ {\sum}^{\Cal N}}  }
 
\centerline{\bf An estimate about multiple stochastic integrals}
\centerline{\bf with respect to a normalized empirical measure}
\smallskip
\centerline{\it P\'eter Major}
\centerline{Alfr\'ed R\'enyi Mathematical Institute of the Hungarian
Academy of Sciences}
\centerline{Budapest, P.O.B. 127 H--1364, Hungary, e-mail:
major\@renyi.hu}
\medskip
{\narrower \noindent {\it Summary:}\/
Let a sequence of iid. random variables $\xi_1,\dots,\xi_n$ be
given on a measurable space $(X,\Cal X)$ with distribution $\mu$
together with a function $f(x_1,\dots,x_k)$ on the product space
$(X^k,\Cal X^k)$. Let  $\mu_n$ denote the empirical measure defined
by these random variables and consider the random integral
$$
J_{n,k}(f)=\frac{n^{k/2}}{k!}\int' f(u_1,\dots,u_k)
(\mu_n(\,du_1)-\mu(du_1))\dots(\mu_n(\,du_k)-\mu(\,du_k)),
$$
where prime means that the diagonals are omitted from the
domain of integration. In this work a good bound is given on the
probability $P(|J_{n,k}(f)|>x)$ for all $x>0$. It is
similar to the bound in the analogous problem we obtain by
considering the Gaussian (multiple) Wiener--It\^o integral of the
function $f$. The proof is based on an adaptation of some methods of
the theory of Wiener--It\^o integrals. In particular, a sort of
diagram formula is proved for the random integrals $J_{n,k}(f)$
together with some of its important properties, a result which may
be interesting in itself. The relation of this paper to some
results about $U$-statistics is also discussed. \par}
 
\beginsection 1. Introduction, formulation of the main results
 
The investigation of the questions discussed in this paper was
motivated by some problems of non-parametric maximum likelihood
estimates. (See~$[5]$,~$[6]$.) The following problem is
studied here: Let $\mu$ be a probability measure on a measurable
space $(X,\Cal X)$, and take a sequence $\xi_1,\dots,\xi_n$ of
independent, identically distributed $(X,\Cal X)$ valued random
variables with distribution~$\mu$. Let us introduce the empirical
measure $\mu_n$,
$$
\mu_n(A)=\dfrac1n\#\{j\:\,\xi_j\in A,\;1\le j\le n\},\quad A\in\Cal X
$$
of this sample $\xi_1,\dots,\xi_n$, and given a function
$f(x_1,\dots,x_k)$ on the $k$-fold product space $(X^k,\Cal X^k)$
define the integral $J_{n,k}(f)$ of the function $f$ with
respect to the $k$-fold product of the normalized empirical
measure $\mu_n$ by the formula
$$
\align
J_{n,k}(f)&=\dfrac{n^{k/2}}{k!} \int'
f(u_1,\dots,u_k)(\mu_n(\,du_1)-\mu(\,du_1))\dots
(\mu_n(\,du_k)-\mu(\,du_k)),\\
&\qquad\text{where the prime in $\tsize\int'$ means that the
diagonals } u_j=u_l,\; 1\le j<l\le k,\\
&\qquad\text{are omitted from the domain of integration.} \tag1.1
\endalign
$$
Let us consider constants as functions of zero variable, and define
$J_{n,0}(c)=c$ for a constant $c$. Given a function $f(x_1,\dots,x_k)$
of $k$ variables we want to give a good estimate on the probability
$P(|J_{n,k}(f)|>x)$ for all $x\ge0$. Our main result
is the following
\medskip\noindent
{\bf Theorem 1.} {\it Let $f=f(x_1,\dots,x_k)$ be a
measurable function on the space $(X^k,\Cal X^k,\mu^k)$ with some
$k\ge1$ such that
$$
\|f\|_\infty=\supp_{x_j\in X,\;1\le j\le k}|f(x_1,\dots,x_k)|\le 1,
$$
and
$$
\|f\|_2^2=Ef^2(\xi_1,\dots,\xi_k)=\int f^2(x_1,\dots,x_k)
\mu(\,dx_1)\dots\mu(\,dx_k)\le \sigma^2 \tag1.2
$$
with some constant $\sigma>0$. Let us also assume that the measure
$\mu$ on $(X,\Cal X)$ is non-atomic. Then there exist some constants
$C=C_k>0$ and $\alpha=\alpha_k>0$, such that the random integral
$J_{n,k}(f)$ defined by formula (1.1) satisfies the inequality
$$
P\(|J_{n,k}(f)|>x\)\le C \max\(\exp\left\{-\alpha
\(\frac x\sigma\)^{2/k}\right\}, \exp\left\{-\alpha
(nx^2)^{1/(k+1)}\right\} \)  \tag1.3
$$
for all $x>0$. These constants $C=C_k>0$ and $\alpha=\alpha_k>0$
depend only on the parameter~$k$.}\medskip
 
It may be useful to reformulate this result in the following equivalent
form: \medskip\noindent
{\bf Theorem 1$'$.} {\it Under the conditions of Theorem 1
$$
P\(|J_{n,k}(f)|>x\)\le C \exp\left\{-\alpha
\(\frac x\sigma\)^{2/k}\right\} \quad \text{for } x\le
n^{k/2}\sigma^{k+1}
$$
with the number $\sigma$ appearing in (1.2) and some universal
constants $C=C_k>0$, $\alpha=\alpha_k>0$, depending
only on the multiplicity~$k$ of the integral $J_{n,k}(f)$.}
\medskip
Theorem 1 clearly implies Theorem~$1'$, since in the case
$x\le n^{k/2}\sigma^{k+1}$ the first term is larger than the second
one in the maximum at the right-hand side of   formula~(1.3). On
the other hand Theorem~$1'$ implies Theorem~1 also if
$x>n^{k/2}\sigma^{k+1}$, since in this case Theorem~$1'$ can be
applied with $\bar\sigma=\(x n^{-k/2}\)^{1/(k+1)}\ge \sigma$. This
yields that $P\(|J_{n,k}(f)|>x\)\le C \exp\left\{-\alpha
\(\dfrac x{\bar\sigma}\)^{2/k}\right\}=C \exp\left\{-\alpha
(nx^2)^{1/(k+1)}\right\}$ if $x>n^{k/2}\sigma^{k+1}$.
 
As it will be seen later, the expression $\|f\|_2^2$ considered in
formula (1.2) has the same order as the variance of the random
integral $J_{n,k}(f)$. Beside this, if $\eta$ is a random
variable with standard normal distribution, then
$P(|\sigma\eta^k|>x)=P\(|\eta|> \(\dfrac x\sigma\)^{1/k}\)\le\const
\exp\left\{-\(\dfrac x\sigma\)^{2/k}\right\}$, and this inequality is
essentially sharp. Thus the results of Theorem~1 or Theorem~$1'$
state that the tail probability $P(|J_{n,k}(f)|>x)$ of the $k$-fold
random integral $J_{n,k}(f)$ behaves similarly to
$P(|\const\sigma\eta^k|>x)$, where the random variable $\eta$ has
standard normal distribution and $\sigma$ is the variance of
$J_{n,k}(f)$, provided that the level~$x$ we consider is less than
$n^{k/2}\sigma^{k+1}$. It can be shown that such a
condition on the level $x$ is really needed.
 
Actually, I am interested in the following more general problem:
Let a nice class of functions $f\in\Cal F$ be given on the space
$(X^k,\Cal X^k)$, and let us give a good bound on the probability
$P\(\supp_{f\in\Cal F}|J_{n,k}(f)|\ge x\)$. In a subsequent paper
[7] I shall show that for nice classes $\Cal F$ of functions,
for instance if $\Cal F$ is a so-called Vapnik--\v{C}ervonenkis
class of functions bounded by 1, a bound similar to that
in Theorem~1 can be given for this maximum. But to show this
first Theorem~1 of this paper has to be proved.
 
I met such a problem when I tried to apply the method of proof for
the existence of a Gaussian limit of the maximum likelihood estimates
to certain non-parametric problems. (See [5] and [6].) The proof in
the parametric case contains a simple but important linearization
step. In this step the function appearing in the maximum likelihood
equation is replaced by its Taylor expansion around the (unknown)
parameter up to the first term, and it is shown that such a
linearization causes only a negligibly small error. In my attempts to
adapt this argument to the non-parametric case the result of Theorem~1
had to be applied. Actually this estimate was needed only in the case
$k=2$. With its help a bound can be obtained which corresponds
to the estimate of the second coefficient in the Taylor expansion of
the classical maximum likelihood equation. Also the problem discussed
in paper~[7] appears in certain non-parametric estimation problems.
For instance in the estimation of a distribution function with the
help of some observations such a problem appears in a natural way. In
this case we have to bound the difference between the distribution
function and its estimate for all numbers $x$, and this requires a
different linearization for all numbers~$x$. In this problem we
want to estimate the supremum of the error in all points~$x$, and this
requires a bound on the supremum of a class of random integrals
$J_{n,2}(f)$.
 
Earlier I could only prove a weaker version of Theorem~1 in
paper~$[4]$ and applied this result. In that paper I could not
prove a better estimate for random integrals $J_{n,k}(f)$ with a
small variance. I do not know of other papers where the
distribution of the random integral $J_{n,k}(f)$ was investigated
directly. On the other hand, some deep and interesting results were
proved about the tail-behaviour of $U$-statistics, more precisely
about the tail-behaviour of so-called degenerated $U$-statistics with
canonical kernel functions (see~$[1]$,~$[2]$), and they are very
close to our results. It may be worth-while to discuss the relation
between them in more detail. To do this first I recall some notions
about $U$-statistics.
 
Let us consider a function $f=f(x_1,\dots,x_k)$ defined on the
$k$-th power $(X^k,\Cal X^k)$ of a space $(X,\Cal X)$ together with
a sequence of independent and identically distributed random
variables $\xi_1,\xi_2,\dots$ which take their
values on this space $(X,\Cal X)$, and let $\mu$ denote their
distribution. We define with their help the
$U$-statistic~$I_{n,k}(f)$
$$
I_{n,k}(f)=\frac1{k!}\summ\Sb 1\le j_s\le n,\; s=1,\dots, k\\
j_s\neq j_{s'} \text{ if } s\neq s'\endSb
f\(\xi_{j_1},\dots,\xi_{j_k}\).   \tag1.4
$$
(The function $f$ in this formula will also be called a kernel
function in the sequel.)
 
A real valued function $f=f(x_1,\dots,x_k)$ on the $k$-th power
$(X^k,\Cal X^k)$ of a space $(X,\Cal X)$ is called a canonical
kernel function (with respect to a probability measure $\mu$
on the space $(X,\Cal X)$) if
$$
\int f(x_1,\dots,x_{j-1},u,x_{j+1},\dots,x_k)\mu(\,du)=0\quad
\text{for all } 1\le j\le k \text{ \ and \ } x_s\in X,  \; s\neq j.
$$
In an equivalent form this means that if $\xi_j$, $1\le j\le k$, are
independent random variables with distribution $\mu$, then
$Ef(\xi_1,\dots,\xi_k|\xi_1=x_1,\dots,\xi_{j-1}=x_{j-1},
\xi_{j+1}=x_{j+1},\dots,\xi_k=x_k)=0$ for all $1\le j\le k$ and
$x_s\in X$, $s\neq j$.
 
The above definition slightly differs from that given in the
literature, where generally it is assumed that the kernel function in
the $U$-statistic is symmetric. But this difference is not essential.
It is not difficult to see that $I_{n,k}(f)=I_{n,k}(\Sym f)$, and our
definition is equivalent to the usual one in the case of symmetric
kernel functions.
 
The definition of the integral $J_{n,k}(f)$ contains a certain kind
of regularization, since integration is taken with respect to the
$k$-th power of the signed measure $\mu_n-\mu$ in it. On the other
hand, the definition of $U$-statistics contains no such
regularization. This has the consequence that only $U$-statistics
with a canonical kernel function satisfy an estimate similar to
Theorem~1. Theorem~1 has the following consequence.
\medskip\noindent
{\bf Theorem 2.} {\it If $f=f(x_1,\dots,x_k)$ is a canonical
function on the space $(X^k,\Cal X^k,\mu^k)$,
$$
\|f\|_\infty=\supp_{x_j\in X,\;1\le j\le k}|f(x_1,\dots,x_k)|\le 1
$$
and
$$
\|f\|_2^2=Ef^2(\xi_1,\dots,\xi_k)=\int f^2(x_1,\dots,x_k)
\mu(\,dx_1)\dots\mu(\,dx_k)\le \sigma^2,
$$
then there exist some constants $C=C_k>0$ and $\alpha=\alpha_k>0$
such that the $U$-statistic $I_{n,k}(f)$ defined in formula (1.4)
satisfies the inequality
$$
P\(n^{-k/2}|I_{n,k}(f)|>x\)\le C \max\(\exp\left\{-\alpha
\(\frac x\sigma\)^{2/k}\right\}, \exp\left\{-\alpha
(nx^2)^{1/(k+1)}\right\} \) \tag1.5
$$
for all $x>0$, where the constants $C=C_k>0$ and $\alpha=\alpha_k>0$
depend only on the parameter $k$.}\medskip
 
I shall show that Theorem~2 is a simple consequence of Theorem 1. It
is also possible to prove Theorem~2 directly and to deduce Theorem~1
from it. But such a deduction demands more work. A basic result of
the theory of $U$-statistics, the Hoeffding decomposition may help,
since it enables one to represent general $U$-statistics (the
random integral $J_{n,k}(f)$ also can be considered as an
$U$-statistic) as a sum of $U$-statistics of different orders with
canonical kernel functions. But to deduce Theorem~1 from Theorem~2
some additional work is still needed. Here I do not deal with
this problem, but the end of Section~2 in [7] indicates what kind of
calculations are needed to work out the details.
 
Part c) of Theorem~2 in $[1]$ gives an estimate equivalent to
our Theorem~2. It states (with the notations of the present paper)
that under the conditions of Theorem~2
$$
P\(n^{-k/2}|I_{n,k}(f)|>x\)\le
c_1\exp\left\{-\dfrac{c_2x^{2/k}}{\sigma^{2/k}
+\(x^{1/k}n^{-1/2}\)^{2/(k+1)}} \right\} \tag1.6
$$
with some universal constants $c_1$  and $c_2$ depending only on
the parameter $k$. In [1] this estimate is called a new
Bernstein-type inequality. It is formally different from our
Theorem~2, but actually they are equivalent if we do not specify
the universal constants in these estimates. It is simple to see
this if the two cases when the first and the second term is dominating
in the denominator at the right-hand side of the expression in (1.6)
are considered separately.
 
\beginsection 2. The idea of the proof
 
We shall prove Theorem 1 by first giving a good estimate on certain
even moments of the random variable $J_{n,k}(f)$ and then by applying
the Markov inequality $P(|J_{n,k}(f)|>x)\le \dfrac{EJ_{n,k}(f)^{2M}}
{x^{2M}}$. In the proof we need a good bound on the moments
$EJ_{n,k}(f)^{2M}$. We do not have to estimate all even moments, but
a good estimate is needed for a sufficiently large class of moments.
It is enough to bound the moments of the form $M=2^m$, $m=0,1,2,\dots$,
and this is done in Proposition~1. The estimate of Proposition~1 is
contained in formula (3.1). This formula presents an upper
bound for $EJ_{n,k}(f)^{2M}$ in the form of the maximum of
two terms, where the first term has the same order as the $2M$-th
moment of a random variable $\const \sigma\eta^k$ with a random
variable $\eta$ with standard normal distribution,
and the second correction term is needed to make relation (3.1)
valid also for small values of $\sigma$.
 
The main problem of this paper is how to prove Proposition~1. To
overcome this difficulty let us consider an analogous problem when the
moments of a multiple stochastic integral with respect to a Gaussian
process (with independent increments) are estimated and try to adapt
the method applied there. An important result about stochastic
integrals, the diagram formula (see e.g. [3]), enables us to prove a
good estimate for the moments of multiple stochastic integrals with
respect to a Gaussian process. It expresses the product of two
multiple stochastic integrals as the sum of certain multiple
stochastic integrals with different multiplicities. Here also a
constant term may appear which can be considered as a stochastic
integral of zero order. The name diagram formula refers
to the fact that the kernel functions of the stochastic integrals
appearing in this identity are defined with the help of certain
diagrams. The diagram formula enables us to express the $2M$-th
power of a multiple Gaussian stochastic integral as the sum of
certain stochastic integrals and to bound its expected value. I
want to show that an appropriate adaptation of this method can
yield the proof of Proposition~1.
 
To carry out such a program we need a result analogous to the diagram
formula where the product of two multiple stochastic integrals
$J_{n,k_1}(f)$ and $J_{n,k_2}(g)$ with respect to a normalized
empirical measure is expressed as the sum of multiple stochastic
integrals. Such a result is proved in Lemma~2, and this also will be
called the diagram formula. The main difference between Lemma~2
and the corresponding result in the Gaussian case is that now some
new terms appear. The kernel functions of the multiple
stochastic integrals in this new version of the diagram formula will be
defined with the help of an object I shall call coloured diagram.
 
This diagram formula enables us to bound the $4M=2^{m+2}$-th moment
of multiple stochastic integrals if we have a good bound on
their $2M=2^{m+1}$-th moment. This suggests to try to give a
bound on the moments $EJ_{n,k}(f)^{2^m}$ by means of an inductive
procedure with respect the parameter $m$. To do this first we have
to estimate the first and second moment of a multiple
stochastic integral $J_{n,k}(f)$. Such an estimate is given in
Lemmas~1 and~3. They are slightly different from the formulas that
appear in the study of Gaussian multiple stochastic integrals.
The expected value of the stochastic integral $J_{n,k}(f)$ may be
non-zero, and instead of an exact identity I can only give an
upper bound for its second moment. The cause of this difference is
that, unlike in the Gaussian case, the (normalized) empirical measures
of disjoint sets are not independent. But the dependence between
them is very weak, therefore the expected value $EJ_{n,k}(f)$
is so small that its non-zero value causes no problem for us. The
bound we give for this expected value does not depend on the sample
size~$n$, at least if $n\ge\frac k2$. On the other hand, the
estimation we gave in Lemma~1 depends on the multiplicity of the
integral $J_{n,k}(f)$. We need such a bound, since we estimate
the moments of the random integrals $J_{n,k}(f)$ simultaneously for
all~$k=1,2,\dots$.
 
The diagram formula reduces the estimate of the $4M=2^{m+2}$-th
moment of a multiple random integral to the estimate of the
$2M=2^{m+1}$-th moment of a sum of multiple random integrals. If we
do not try to get the best possible constant $\alpha$ in Theorem~1
then it is enough to give a sufficiently good estimate on the $2M$-th
moment of each term in this sum. To get such an estimate the
following analytical problem has to be solved: Given two functions
$f(x_1,\dots,x_{k_1})$ and $g(x_1,\dots,x_{k_2})$ let us express
the product of the random integrals $J_{n,k_1}(f)J_{n,k_2}(g)$ by
means of the diagram formula and give a good estimate on the
$L_2$-norm of the kernel functions appearing at the right-hand side
of this expression. The solution of this problem is a most essential
part of our proof.
 
Such an estimate is given in formulas (3.14) and (3.15) which can be
proved in a natural way by means of the Schwarz inequality. Formula
(3.14) yields an estimate for the $L_2$ norm of general kernel functions
appearing in the diagram formula, and formula (3.15) provides a
better bound in the case of those kernel functions which also appear
in the Gaussian case. But these estimates alone are not sufficient
for our purposes. It turned out that in certain cases the $L_2$ norm
of the kernel functions defined in the diagram formula can be better
bounded. Moreover, to get a sufficiently good estimate for the
moments bounded in Proposition~1 we have to exploit this improvement.
 
This lead us to the introduction of the notion of $(r,\sigma^2)$
dominated functions and the formulation of Lemma~4 which contains the
most important properties of functions with such a property. After
Lemma~4 I also formulated Proposition~2. This is a sharpened version
of Proposition~1 which yields a better estimate for the expected
value $EJ_{n,k}(f)^{2M}$ if the variance $\sigma^2$ of $J_{n,k}(f)$
is small and $f$ is an $(r,\sigma^2)$ dominated function with some
large~$r$. Proposition~2 was introduced not because the estimate
of Proposition~1 was not good enough for us, but because the inductive
procedure we apply in the proof works only for Proposition~2. To
carry out this procedure we must have a sufficiently good estimate
for each term appearing in the intermediate steps, and only
Proposition~2 provides a good estimate for them.
 
To prove Proposition~2 still we need a technical result formulated
in Lemma~5. By working out the arguments sketched above the expected
value $EJ_{n,k}(f)^{2M}$ can be estimated by a sum with terms
possessing a good upper bound. We have to show that these bounds also
yield a sufficiently good estimate for the sum itself. This is not a
trivial problem, and Lemma~5 helps us to overcome this difficulty.
 
This paper consists of five sections. In Sections~1 and 2 the main
results of the paper were presented and the method of the proof was
explained. The lemmas and propositions needed in the proof are
formulated Section~3, and also Theorems~1 and~2 are proved there.
Lemmas 1---4 are proved in Section~4. They describe some important
properties of the random integrals $J_{n,k}(f)$ defined in formula
(1.1). I also put Lemma 4 to Section~4, although it contains a result
about a deterministic transform of some functions, so formally this
result is not connected to the properties of the random integrals
$J_{n,k}$. Nevertheless, it can be put among the important results
about such random integrals, because it plays an essential role in the
estimation of the high moments of the terms appearing in the diagram
formula. Finally, Section~5 contains the proof of Proposition~2
together with the rather technical Lemma~5 needed in its proof.
Theorem~1 was proved as a consequence of this Proposition.
 
\beginsection
3. The proof of Theorems 1 and 2 with the help of some lemmas
 
Theorem~1 will be proved with the help of the following result.
\medskip\noindent
{\bf Proposition 1.} {\it If $f=f(x_1,\dots,x_k)$ is a function of
$k$ variables,
$$
\|f\|_\infty=\supp_{x_j\in X,\;1\le j\le k}|f(x_1,\dots,x_k)|\le 1,
$$
and
$$
\|f\|_2^2=\int f^2(x_1,\dots,x_k) \mu(\,dx_1)\dots\mu(\,dx_k)\le
\sigma^2,
$$
the measure $\mu$ is non-atomic, then there exist some constants
$C_k$ such that for all numbers of the form $M=2^m$,
where $m$ is a non-negative integer, and $n=1,2,\dots$
$$
E J_{n,k}(f)^{2M} \le \(C_k\sigma^2 M^k\)^M\cdot \max\(1,
\(\frac M{n\sigma^2}\)^M \). \tag3.1
$$
The constant $C_k$ depends only on the number of variables $k$ in the
function $f$. } \medskip
Proposition~1 states that typically the $2M$-th moment of the random
variable $J_{n,k}(f)$ behaves like the $2M$-th moment of $\const
\sigma\eta^k$, where $\eta$ is a random variable with standard
normal distribution. But such an estimate holds only if $n\sigma^2$
is not too small. Formula~3.1 also contains a constant $C_k$ not
given explicitly. The parameter $\alpha$ in the exponent of the
estimate in Theorem~1 depends on its value.
\medskip
To prove Proposition 1 we give an upper bound for the first and
second moment of a random integral $J_{n,k}(f)$ if $f$ is a bounded
function with $k$ variables and prove some kind of diagram formula
which expresses the random integral $J_{n,k}(f)^2$ or more generally
the product of random integrals $J_{n,k_1}(f)J_{n,k_2}(g)$ with
bounded functions $f$ and $g$ of $k_1$ and $k_2$ arguments as the
sum of appropriate random integrals $J_{n,\bar k}(\cdot)$ of some
functions which have $0\le \bar k\le k_1+k_2$ arguments.
 
To express a product $J_{n,k_1}(f)J_{n,k_2}(g)$ of random integrals
with functions $f$ and $g$ of $k_1$ and $k_2$ arguments respectively
in the desired form we introduce some diagrams and coloured diagrams
and define some functions with their help.
 
Given two positive integers $k_1>0$ and $k_2>0$ let us consider a set
$$
\Cal N=\{(j_1,j_1'),\dots,(j_l,j_l')\}
$$
of pairs of integers such that $1\le j_s\le k_1$, $k_1+1\le j_s'\le
k_1+k_2$, $1\le s\le l$, and $j_s\neq j_{s'}$, $j'_s\neq j'_{s'}$ if
$s\neq s'$. Let us assume that $1\le j_1<j_2<\cdots<j_l\le k_1$ to
introduce an explicit ordering among the elements of the set $\Cal N$.
We define the diagram $B(\Cal N)=B(\Cal N,k_1,k_2)$ in the following
way: The diagram $B(\Cal N)$ has two rows, the first row consists
of the vertices $1,\dots,k_1$ and the second row of the vertices
$k_1+1,\dots,k_1+k_2$. The diagram $B(\Cal N)$ contains the edges
connecting the vertex $j_s$ from the first row with the vertex $j'_s$
in the second row, $(j_s,j_s')\in \Cal N$, $1\le s\le l$. Given
a diagram $B(\Cal N)$ and a subset $\Cal N_1\subset \Cal N$ we
define the coloured diagram $B(\Cal N,\Cal N_1)=B(\Cal N,\Cal N_1,
k_1,k_2)$ as the diagram $B(\Cal N)$ whose edges are coloured in the
following way: An edge $(j_u,j'_u)$ has colour 1 if $(j_u,j'_u)\in
\Cal N_1$, and colour $-1$ if $(j_u,j'_u)\in\Cal N\setminus \Cal N_1$.
 
Let us also define the following operators on the space of functions.
Given an integrable function $h(x_1,\dots,x_s)$ on the product space
$(X^s,\Cal X^s,\mu^s)$ with some $s=1,2,\dots$, define the operators
$P_j=P_{j,s}$, $1\le j\le s$, as
$$
P_jh(x_1,\dots,x_{j-1},x_{j+1},\dots,x_s)=\int
h(x_1,\dots,x_{j-1},u,x_{j+1},\dots,x_s)\mu(\,du). \tag3.2
$$
Let us fix a pair of integers $(j,j')$, $1\le j\le k_1$,
$k_1+1\le j'\le k_1+k_2$ and a function $h(x_{s(1)},\dots,x_{s(r)})$
on an appropriate product of the space $(X,\Cal X,\mu)$ whose
variables are indexed by some numbers $1\le s(1)<s(2)<\cdots<s(r)\le
k_1+k_2$, and $\{j,j'\}\subset\{s(1),\dots,s(r)\}$. In this case we
define the function $R_{j,j'}h$ as
$$
\aligned
&R_{j,j'}h(x_{s(1)},\dots,x_{s(p-1)},x_{s(p)},x_{s(p+1)},\dots,
x_{s(q-1)},x_{s(q)},x_{s(q+1)},\dots,x_{s(r)}) \\
&\qquad =h(x_{s(1)},\dots,x_{s(p-1)},x_{s(p)},x_{s(p+1)},\dots,
x_{s(q-1)},x_{s(p)},x_{s(q+1)},\dots,x_{s(r)}) \\
&\qquad\qquad \text{if } s(p)=j,\text{ and } s(q)=j'
\endaligned\tag3.3
$$
In words: The application of the operator $R_{j,j'}$ means that
the argument $x_{j'}$ is replaced by the argument $x_j$.
 
Given two functions $f(x_1,\dots,x_{k_1})$ and
$g(x_1,\dots,x_{k_2})$ on the spaces $(X^{k_1},\Cal X^{k_1})$
and $(X^{k_2},\Cal X^{k_2})$ respectively define the function
$$
f\circ g(x_1,\dots,x_{k_1+k_2})=
f(x_1,\dots,x_{k_1})g(x_{k_1+1},\dots,x_{k_1+k_2}) \tag3.4
$$
and given a diagram $B(\Cal N)=B(\Cal N,k_1,k_2)$ with edges
$\Cal N=\{(j_1,j_1'),\dots,(j_l,j_l')\}$
put
$$
\overline {f\circ g}_{B(\Cal N)}(x_{s(1)},\dots,x_{s(k_1+k_2-l)})
=\prod_{t=1}^l R_{j_t,j'_t} f\circ g(x_1,\dots,x_{k_1+k_2})
\tag3.5
$$
with $s(u)=u$ if $1\le u\le k_1$ and to define $s(u)$ in
the case $k_1+1\le u\le k_1+k_2-l$ let us list the set
$\{k_1+1,\dots,k_1+k_2\}\setminus \{j_1',\dots,j_l'\}
=\{v(1),\dots,v(k_2-l)\}$, $k_1+1\le v(1)<\cdots <v(k_2-l)\le k_1+k_2$
in increasing order and put $s(u)=v(u-k_1)$. In words, the following
construction was made. Take the product of the functions $f$ and $g$,
and replace the argument $x_{j_s'}$ by $x_{j_s}$ for all edges
$(j_s,j_{s'})$ of the diagram $B(\Cal N)$, and then the indices of
the variables are not changed. The numbers $s(u)$ were introduced to
list the indices of the variables in an increasing order.
 
Given a coloured diagram
$B(\Cal N,\Cal N_1)=B(\Cal N,\Cal N_1,k_1,k_2)$ with
$$
\Cal N=\{(j_1,j_1'),\dots,(j_l,j_l')\} \text{ and }
\Cal N_1=\{(j_{u(1)},j_{u(1)}'),\dots,(j_{u(p)},j_{u(p)}')\},
$$
where $\{u(1),\dots,u(p)\}\subset\{1,\dots,l\}$ define
$$
\overline {f\circ g}_{B(\Cal N,\Cal N_1)}
(x_{s(1)},\dots,x_{s(k_1+k_2-l-p)})
=\prod_{v=1}^p P_{j_{u(v)}} \prod_{t=1}^l R_{j_t,j'_t}
f\circ g(x_1,\dots,x_{k_1+k_2}) \tag3.6
$$
with the help of the operators $P_j$ and $R_{j,j'}$ defined in
(3.2) and (3.3), where the operators at the right-hand side are applied
from right to left order and then the indices $u(\cdot)$ at the
left-hand side are defined in the following way: Let us list the sets
$\{w(1),\dots,w(k_1-p)\}=\{1,\dots,k_1\}\setminus
\{j_{u(1)},\dots,j_{u(p)}\}$, $1\le w(1)<\cdots<w(k_1-p)$, and
$\{v(1),\dots,v(k_2-l)\}= \{k_1+1,\dots,k_1+k_2\}\setminus
\{j_1',\dots,j_l'\}$, $k_1+1\le v(1)<\cdots <v(k_2-l)\le k_1+k_2$,
in increasing order. Then we define $s(u)=w(u)$ if $1\le u\le k_1-p$,
and $s(u)=v(u-k_1+p)$ if $k_1-p+1\le u\le k_1+k_2-l-p$. In words, take
the function $\overline {f\circ g}_{B(\Cal N)}$ and apply the
projection operator $P_u$ for those indices $u$ which are end-points
of edges with colour 1 in the first row of the coloured diagram
$B(\Cal N,\Cal N_1)$  i.e. which are end-points  of edges in $\Cal
N_1$ in the first row. After this procedure these variables $x_u$
disappear. The remaining variables  are preserved together with
their indices. The numbers $s(u)$ were
introduced again to list these variables in increasing order.
 
To formulate the diagram formula we define the function
$f\circ g_{B(\Cal N,\Cal N_1)}$ with the help of the function
$\overline{f\circ g}_{B(\Cal N,\Cal N_1)}$ in the following way.
$$
{f\circ g}_{B(\Cal N,\Cal N_1)} (x_1,\dots,x_{k_1+k_2-l-p})=
\overline {f\circ g}_{B(\Cal N,\Cal N_1)}
(x_{s(1)},\dots,x_{s(k_1+k_2-l-p)}). \tag3.7
$$
Here we reindexed the variables of $\overline {f\circ g}$ and listed
them in an increasing order.
 
Now I formulate some lemmas needed in the proof of Proposition~1.
\medskip\noindent
{\bf Lemma~1.} {\it Let $f(x_1,\dots,x_k)$ be an integrable
function on a space $(X^k,\Cal X^k,\mu^k)$ where $\mu^k$ is the
$k$-th power of the probability measure $\mu$ appearing in formula
(3.1). Let the measure $\mu$ be non-atomic. Then
$$
EJ_{n,k}(f)=B_{n,k}\int f(u_1,\dots,u_k)\mu(\,du_1)\dots\mu(\,du_k)
\tag3.8
$$
with
$$
B_{n,k}=\frac1{k! n^{k/2}}\sum_{s=1}^k (-1)^{k-s}\binom ns s!
\sum\Sb (r_1,\dots,r_s)\\ r_j\text{ is integer } r_j\ge1, \;1\le j\le
s,\\ s_1+\cdots+r_s=k\endSb
(r_1-1)\cdots(r_s-1)B(r_1,\dots,r_s),
\tag3.9 $$
where $B(r_1,\dots,r_s)$ equals the number of partitions of the set
$\{1,\dots,k\}$ to disjoint sets with cardinalities $r_1,\dots,r_s$.
The above constants $B_{n,k}$ satisfy the estimate
$$
|B_{n,k}|\le \dfrac{C^k}{k^{k/2}} \tag3.10
$$
with some universal constant $C>0$ for all $n\ge\frac k2$.}
\medskip\noindent
{\bf Lemma 2. (Diagram formula)} {\it Let $f=f(x_1,\dots,x_{k_1})$
and $g(x_1,\dots,x_{k_2})$ be two square integrable
functions on the product spaces $(X^{k_1},\Cal X^{k_1},\mu^{k_1})$
and $(X^{k_2},\Cal X^{k_2},\mu^{k_2})$ respectively, and assume that
the measure $\mu$ is non-atomic. Then the identity
$$
J_{n,k_1}(f)J_{n,k_2}(g)=\!\!  \sum_{l=0}^{\min(k_1,k_2)} \!
\sum_{p=0}^l \frac{(k_1+k_2-l-p)!}{(k_1-l)!(k_2-l)!(l-p)!p!}
n^{-(l-p)/2}
J_{n,k_1+k_2-l-p}(f\circ g_{l,p}). \tag3.11
$$
holds with the functions $f\circ g_{l,p}=f\circ g_{l,p,k_1,k_2}$ of
$k_1+k_2-l-p$ variables, $0\le p\le l\le\min(k_1,k_2)$, defined by
the formula
$$
\aligned
f\circ g_{l,p}(x_1,\dots,x_{k_1+k_2-l-p}) &=
\frac{(k_1-l)!(k_2-l)!(l-p)!p!}{k_1! k_2!} \\
&\quad \sum_{B(\Cal N,\Cal N_1)\in\Cal B(l,p)}
f\circ g_{B(\Cal N,\Cal N_1)}(x_1,\dots,x_{k_1+k_2-l-p})
\endaligned \tag3.12
$$
where the functions $f\circ g_{B(\Cal N,\Cal N_1)}$ are defined by
formulas (3.2), (3.3), (3.6) and (3.7), and the set $\Cal B(l,p)=
\Cal B(l,p,k_1,k_2)$ consists of those coloured diagrams $B(\Cal
N,\Cal N_1)=B(\Cal N,\Cal N_1,k_1,k_2)$ with $k_1$ vertices in the first
and $k_2$ vertices in the second row for which $|\Cal N|=l$,
and $|\Cal N_1|=p$. (We apply here and in the sequel the convention
that $|A|$ denotes the cardinality of a set $A$.)}
\medskip\noindent
Let us remark that the norming constant
$\dfrac{(k_1-l)!(k_2-l)!(l-p)!p!}{k_1!k_2!}$ in (3.12) is the
reciproque of the cardinality of the set $\Cal B(l,p,k_1,k_2)$.
\medskip\noindent
{\bf Lemma 3.} {\it Let $f(x_1,\dots,x_k)$ be a square integrable
function on the space $(X^k,\Cal X^k)$, and let $\mu$ be a non-atomic
measure on $(X,\Cal X)$. Then for all $n\ge k$
$$
EJ_{n,k}(f)^2\le \dfrac {C^k}{k^k} \|f\|_2^2 \tag3.13
$$
with some universal constant $C>0$.}
\medskip
The diagram formula formulated in Lemma 2 enables us to estimate
$EJ_{n,k}(f)^{2^{m+1}}$ if we have an estimate for the expectation
$EJ_{n,\bar k}(g)^{2^m}$ for all bounded functions $g$ (where also
the number of arguments in the function $g$ may be arbitrary). This
may enable us to prove Proposition~1 by induction. But to prove a
sufficiently good estimate, first we have to prove a good bound on
the $L_2$ norm of functions of the form
$f\circ f_{B(\Cal N,\Cal N_1)}(x_1,\dots,x_{2k-p})$ with some coloured
diagram $B(\Cal N,\Cal N_1)$. We will show that
$$
\|f\circ g_{B(\Cal N,\Cal N_1)}\|_2^2\le \|f\|_2\|g\|_2 \quad \text
{for any diagram } B(\Cal N,\Cal N_1)
\text{ if } \|f\|_\infty\le1 \text{ and } \|g\|_\infty\le1.
\tag3.14
$$
Moreover, in the special case when $|\Cal N|=|\Cal N_1|$
this estimate can be improved. Namely, we shall show that
$$
\|f\circ g_{B(\Cal N,\Cal N_1)}\|_2^2\le
\|f\|_2^2 \|g\|_2^2\quad \text{if } |\Cal N|=|\Cal N_1|. \tag3.15
$$
Let us call a diagram $B(\Cal N,\Cal N_1)$ Gaussian if $\Cal N=\Cal
N_1$, and non-Gaussian in the other case. The introduction of such
a terminology is natural, since in the diagram formula for multiple
stochastic integrals with respect to Gaussian processes those
diagrams appear (with the same integrals attached to them) which we
called Gaussian. The most difficult part of our investigation is to
check the contribution of the integrals corresponding to
non-Gaussian diagrams.
 
To prove formula (3.14) let us take two functions
$f(x_1,\dots,x_{k_1})$ and $g(x_1,\dots,x_{k_2})$ whose supremum
norm is bounded by~1 and consider the function $\overline{f\circ
g}_{B(\Cal N)}$ with some diagram $\Cal N$ defined by formulas (3.3),
(3.4) and (3.5). The Schwarz inequality yields that
$$
\align
\left\|\overline {f\circ g}_{B(\Cal N)}\right\|_2^2 &\le \(\int
f^4(x_1,\dots,x_k)\mu(\,dx_1) \dots \mu(\,dx_{k_1})\)^{1/2}\\
&\qquad \(\int g^4(x_1,\dots,x_{k_2})\mu(\,dx_1)\dots
\mu(\,dx_{k_2})\)^{1/2}\\
&\le \(\int f^2(x_1,\dots,x_k)\mu(\,dx_1) \dots
\mu(\,dx_{k_1})\)^{1/2}\\
&\qquad\(\int g^2(x_1,\dots,x_{k_2})\mu(\,dx_1)\dots
\mu(\,dx_{k_2})\)^{1/2}
=\|f\|_2 \|g\|_2
\endalign
$$
if $\|f\|_\infty\le1$ and $\|g\|_\infty\le1$.
 
The operators $P_j$ are contractions in $L_2$ norm. Hence the above
estimate together with formulas (3.5), (3.6) and (3.7) imply that
$$
\left\|f\circ g_{B(\Cal N,\Cal N_1)}\right\|_2^2=
\left\|\overline{f\circ g}_{B(\Cal N,\Cal N_1)}\right\|_2^2
\le\left\|\overline{f\circ g}_{B(\Cal N)}\right\|_2^2
\le\|f\|_2 \|g\|_2,
$$
and this is relation (3.14)
 
To prove formula (3.15) let us first give a pointwise estimate on
$\overline{f\circ g}_{B(\Cal N,\Cal N)}^2$ with the help of the
Schwarz inequality. We get by using formulas (3.2),
(3.3), (3.4) and (3.6) that
$$
\aligned
\overline {f\circ g}_{B(\Cal N,\Cal N)}^2
(x_{s(1)},\dots,x_{s(k_1+k_2-l-p)})
&\le\int f^2(x_1,\dots,x_{k_1})\prod_{s=1}^l\mu(\,dx_{j_s}) \\
&\qquad \int g^2(x_{k_1+1},\dots,x_{k_1+k_2})\prod_{s=1}^l
\mu(\,dx_{j_s'}), \endaligned \tag3.16
$$
where $\Cal N=\{(j_1,j_1'),\dots,(j_l,j_l')\}$. Beside this, the
arguments of the function at the left-hand side of (3.16) are the
elements of the set
$\{s(1),\dots,s(k_1+k_2-l-p)\}=\{1,\dots,k_1+k_2\}
\setminus\{j_1,\dots,j_l,j_1',\dots,j_l'\}$.
The important point is that the estimate at the right-hand side of
(3.16) is the product of two functions with different
arguments. We get by integrating both sides of inequality (3.16)
that
$$
\left\|\overline {f\circ g}_{B(\Cal N,\Cal N)}\right\|^2\le
\|f\|_2^2\|g\|_2^2.
$$
This inequality together with the identity
$\|\overline {f\circ g}_{B(\Cal N,\Cal N)}\|^2_2=
\|f\circ g_{B(\Cal N,\Cal N)}\|^2_2$ imply formula (3.15).
 
Lemma 3 gives a bound on $EJ_{n,k}(f)^2$ for a function $f$ with the
help of $\|f\|_2^2$, and the application of Lemma~2 with the choice
$f=g$ enables us to bound $EJ_{n,k}(f)^{4M}$ if we have a bound for
the expressions of the form $EJ_{n,\bar k}(g)^{2M}$. In such a way
we can prove good bounds on the moments $EJ_{n,k}(f)^{2^m}$ by
induction. To carry out such a program we need a good bound on the
$L_2$ norm of the terms at the right-hand side (3.12) in the special
case $f=g$. Formulas (3.14) and (3.15) supply such estimates, but
they are not always good enough for our purposes.
 
These estimates are sufficient to bound $EJ_{n,k}(f)^{2M}$ for
$M=1,2$, but do not suffice already in the case $M=4$. Indeed, to
give a good bound on $EJ_{n,k}(f)^8$ with the help of Lemma~2 we
also have to give a good bound on $EJ_{n,2k}(g)^4$ with
$g(x_1,\dots,x_{2k})=f(x_1,\dots,x_k)f(x_{k+1},\dots,x_{2k})$. (This
term is $f\circ f_{B(\Cal N,\Cal N_1)}$ with $\Cal N=\Cal
N_1=\emptyset$.) To estimate this expression (again with the help of
Lemma 2) we also have to bound $EJ_{n,4k-1}(g_1)^2$ with the function
$g_1(x_1,\dots,x_{4k-1})=g(x_1,\dots,x_{2k})
g(x_1,x_{2k+2},\dots,x_{4k})$ (here we consider the diagram $B(\Cal
N,\Cal N_1)$ with $\Cal N=\{(1,2k+1)\}$ and $\Cal N_1=\emptyset$),
and to estimate $EJ_{n,4k-1}(g_1)^2$ we  have to bound $L_2$ norm of
the function~$g_1$.
 
Formulas (3.14) and (3.15) imply that $\|g\|_2^2\le\sigma^4$
and $\|g_1\|_2^2\le\sigma^4$ if $\|f\|_2^2\le\sigma^2$. But a more
careful analysis shows that the function $g_1$ also satisfies the
sharper estimate $\|g_1\|_2^2\le \sigma^6$. Indeed, the function $g$
is the product of two functions with different arguments, and the
$L_2$ norm of both terms is bounded by $\sigma$. Further
considerations show that the function $g_1$ can be written as the
product of three terms whose $L_2$ norm are bounded by $\sigma$,
and this implies the above statement.
 
If $\|f\|_2^2\le\sigma^2$ with a small number $\sigma$, then to
get a sufficiently good estimate for $EJ_{n,k}(f)^{2M}$ for
$M\ge 4$ with the help of Lemma~2 a sufficiently good bound has to
be proved about $L_2$ norm of certain functions. The estimates
(3.14) and (3.15) are not always sufficient for this goal. In the
previous example we could get some improvement by exploiting the
special structure of certain functions. It worked, because the
function whose $L_2$ norm we had to bound could be represented as
the product of functions with different arguments and small $L_2$
norms. An appropriate refinement of this argument would also supply
a useful estimate for the $L_2$ norm of such functions which can be
bounded by the product of functions with different arguments and
small $L_2$ norm. To get sufficiently good estimates which lead to
the proof of Propitiation~1 we have to exploit the improvements
following from the above sketched observations. Such considerations
lead to the introduction of the following definition.
\medskip\noindent
{\bf The definition of $(r,\sigma^2)$ dominated functions.} {\it Let
us consider a function $f=f(x_{u_1},\dots,x_{u_k})$, $x_{u_j}\in X$,
$1\le j\le k$, of $k$ variables. We say that this function is
$(r,\sigma^2)$ dominated with $1\le r<\infty$, $0<\sigma\le1$ if the
set $\{u_1,\dots,u_k\}$ can be partitioned to $r$ disjoint sets
$B_l=\{{j_1}^{(l)},\dots,{j_{s_l}}^{(l)}\}$, $1\le l\le r$,
$\bigcupp_{l=1}^r B_l=\{u_1,\dots,u_k\}$, and for all $1\le l\le r$
such functions $h_l=h_l(x_{j_1^{(l)}},\dots,x_{j_{s_l}^{(l)}})$
can be defined whose arguments are indexed by the elements of the
set $B_l$, and which satisfy the inequalities
$$
|f(x_{u_1},\dots,x_{u_k})|\le \prod_{l=1}^r
h_l(x_{j_1^{(l)}},\dots,x_{j_{s_l}^{(l)}})
\quad \text{for all } (x_{u_1},\dots,x_{u_k})
$$
and $\|h_l\|_\infty\le 1$, $\|h_l\|_2^2\le \sigma^2$, $1\le l\le r$.
 
In particular, we allow that some of the sets $B_l$ be empty in the
above partition. If the set $B_l$ is empty, then the function $h_l$
has to be constant, and $0\le h_l\le \sigma$.} \medskip
Let us remark that a function $f$ is $(1,\sigma^2)$ dominated if and
only if $\|f\|_\infty\le1$, and $\|f\|_2\le \sigma$.
 
The following Lemma~4 which uses the above definition plays a most
important role in the proof of Proposition~1. We formulate it in
a form more general than we need. We shall consider two functions
$f$ and $g$, although we shall apply this result only in the special
case when $f=g$. \medskip\noindent
{\bf Lemma 4.} {\it Let $f$ be a function on the product space
$(X^{k_1},\Cal X^{k_1},\mu^{k_1})$ and $g$ a function on the product
space $(X^{k_2},\Cal X^{k_2},\mu^{k_2})$ with some integers $k_1\ge1$,
$k_2\ge1$. Let us take a coloured diagram $B(\Cal N,\Cal N_1)=B(\Cal
N,\Cal N_1,k_1,k_2)$ with $|\Cal N|=l$ and $|\Cal N_1|=p$ and
consider the function $f\circ g_{B(\Cal N,\Cal N_1)}$ defined
by formulas (3.2), (3.3), (3.6) and (3.7). Let $f$ be an
$(r_1,\sigma^2)$ and $g$ an $(r_2,\sigma^2)$ dominated function.
If the relation $r_1+r_2\ge l-p+1$ holds, then the function
$f\circ g_{B(\Cal N,\Cal N_1)}$ is $(r_1+r_2-(l-p),
\sigma^2)$ dominated. Beside this, the function $f\circ g_{B(\Cal
N,\Cal N_1)}$ is $(1,\sigma^{(r_1+r_2)/2})$ dominated for all
$0\le p\le l\le\min(k_1,k_2)$.}
\medskip
 
It is clear that if $f$ is an $(r_1,\sigma)$ and $g$ is an
$r_2,\sigma)$ dominated function, then $f\circ g$ is an $r_1+r_2,
\sigma$ dominated function. Lemma~4 essentially states that if we
calculate the function $f\circ g_{B(\Cal N,\Cal N_1)}$ with the
help of the function $f\circ g$, then the edges of $\Cal N_1$
(with colour 1) do not diminish and the edges in $\Cal N\setminus
\Cal N_1$ diminish only with 1 the first coordinate in the
$(r_1+r_2,\sigma^2)$  dominance property.
 
With the help of Lemmas~1---4 the following Proposition~2 can be
proved which yields a more detailed description on the moments
of the random integrals $J_{n,k}(f)$ than Proposition~1.
\medskip\noindent
{\bf Proposition~2.} {\it Let $f=f(x_1,\dots,x_k)$ be a measurable
$\(r,\sigma^{2/r}\)$ dominated function of $k$ variables, $k\ge1$,
with some $r\ge1$. Put $M=2^m$. Let us assume that the measure $\mu$
appearing in formula (3.1) is non-atomic. Then for all $m=0,1,\dots$
such that $kM\le n$
$$
EJ_{n,k}(f)^{2M} \le \(\frac{C(k,m) M^k \sigma^2}{k^k}\)^M
\cdot \max \(1,\(\frac{kM}{n\sigma^{2/r}} \)^{M\min(k,r)}\)
\tag3.17
$$
with some constants $C(k,m)\ge1$ which do not depend on the
parameters $n$ and $\sigma$, and they satisfy the inequality
$\supp_{0\le m<\infty} C(k,m)\le C_k<\infty$ for all $k=1,2,\dots$.}
\medskip\noindent
{\it Remark:}\/ The condition $kM\le n$ in Proposition 2 is of
technical character which probably can be omitted. We imposed this
condition to avoid the investigation of such degenerate random
integrals $J_{n,k}(f)$ where the number of variables of the function
$f$ is larger than the sample size $n$. To avoid such a situation
during the inductive procedure leading to the proof of Proposition~2
we had to impose the condition~$kM\le n$. This restriction does not
cause a serious problem in the proof of Theorem~1.
\medskip
Proposition 1 follows from Proposition~2 if we apply it for
$\(r,\sigma^{2/r}\)$ dominated functions $f=f(x_1,\dots,x_k)$
with $r=1$ and choose $C=C_k=\dfrac1{k^{k-1}}\supp_{0\le m\le \infty}
C(k,m)$. Actually Proposition 2 implies Proposition~1 only under the
additional condition $kM\le n$. But Proposition~1 can be proved in
the case $n\le kM$ with the help of the following very rude estimate:
Since the measure $\mu_n-\mu$ is the difference of two probability
measures and $\|f\|_\infty\le 1$, hence $\|J_{n,k}f\|_\infty\le 2^k
n^{k/2}$ and $EJ_{n,k}(f)^{2M}\le 2^{2kM} n^{kM}\le 2^{2kM}(kM)^{kM}
\le n^{-M} 2^{2kM}(kM)^{(k+1)M}=\bar C_k^M n^{-M}M^{(k+1)M}$ with
$\bar C_k=4^k k^{k+1}$ if $n\le kM$. This estimate means that the
estimate (3.1) remains valid with a possibly new constant $C_k$ if
we also allow the case $n\le kM$. (Here we use the second term in
the maximum at the right-hand side of (3.1).) Hence it is enough to
prove Proposition~2.
 
Proposition 2 will be proved by induction with respect to $m$. For
$m=0$ Lemma~3 implies the result we need. If we know the bound (3.17)
for some $m$ with a constant $C(k,m)$, then this bound together
with Lemma~4 enable us to give a good bound on the $2^m$-th
moment of the expressions $J_{n,2k-(l-p)}(f\circ f_{l,p})$ with the
functions $f\circ f_{l,p}$ appearing in formula (3.12). Then we can
prove formula (3.17) for $m+1$ with an appropriate constant $C(k,m+1)$
by means of the diagram formula presented in Lemma~2. The calculations
made during this proof show that constant $C(k,m+1)$ appearing in this
estimate depends only on the parameters $k$ and $m$, i.e.\ it does not
depend on the sample size $n$ or the variance $\sigma^2$.
Moreover, an important but rather technical result formulated below
in Lemma~5 will enable us to show that the norming constants $C(k,m)$
can be chosen relatively small. By exploiting the existence of some
constants $\bar C(k,m)$ with the properties formulated in Lemma~5 we
can show that in Proposition~2 the constants $C(k,m)$ can be
chosen in such a way that they also satisfy the inequality
$\supp_{0\le m\le \infty}C(k,m)<\infty$. \medskip\noindent
{\bf Lemma 5.} {\it There exists a set of positive real numbers
$\bar C(k,m)$, $\bar C(k,0)=1$, $m=0,1,\dots$, $k=0,1,\dots$, which
satisfy the inequality
$$
\aligned
\bar C(k,m+1)^2&\ge
2^{2l(4-m)} \frac{(2k)^{2k-l+p}(2k-l-p)^{3l-p-2k}}
{(2l)^{2l}} \bar C(2k-l-p,m)\\
&\qquad\text{for all } k=0,1,2,\dots,\; m=0,1,2,\dots \text { and }
0\le p\le l\le k,
\endaligned \tag3.18
$$
and
$$
\sup_{m\ge0}\bar C(k,m)\le
D(k)<\infty\quad\text {for all }k=0,1,2,\dots \tag3.19
$$
with some constants $D(k)>0$, $k=0,1,2,\dots$. The numbers $\bar
C(k,m)$ can be given in the form $\bar C(k,m)=\bar
C(m)^k$, where $\bar C(m)$, $m=0,1,2,\dots$, is a monotone increasing
sequence such that $\bar C(0)=1$, and $\limm_{m\to\infty}\bar
C(m)<\infty$.}\medskip\noindent
 
We finish Section~2 by proving Theorems~1 and~2 with the help of
Proposition~1. \medskip\noindent
{\it The proof of Theorem 1.}
By Proposition 1 and the Markov inequality
$$
P\(|J_{n,k}(f)|>x\)\le
\(\frac{C_k\sigma^2 M^k}{x^2}\)^M\cdot \max\(1,
\(\frac M{n\sigma^2}\)^M \) \tag3.20
$$
for all positive integers $M$ of the form $M=2^m$. Let us assume
for a while that the number $x$ satisfies the inequalities $\dfrac
x\sigma\ge 2^{(k+1)/2}C_k^{1/2}$ and $nx^2\ge 2^{k+3}C_k$. Then the
number $M=2^m$ can be chosen in such a way that the inequality
$\dfrac1{2^{k+1}}\le \dfrac{C_k\sigma^2 M^k}{x^2}<\dfrac12$, i.e.\
$\dfrac1{(2^{k+1}C_k)^{1/k}}\(\dfrac x\sigma\)^{2/k}\le M<
\dfrac1{(2C_k)^{1/k}}\(\dfrac x\sigma\)^{2/k}$ holds. Let us choose
such a number $M$ if $\dfrac1{(2C_k)^{1/k}}\(\dfrac
x\sigma\)^{2/k}\le n\sigma^2$, i.e.\ $x\le
\sqrt{2C_k}n^{k/2}\sigma^{k+1}$. In this case $M\le n\sigma^2$,
the first term is dominating in the maximum at the right-hand side
of formula (3.20), and $P\(|J_{n,k}(f)|>x\)\le 2^{-M}\le
\exp\left\{-\alpha\(\dfrac x\sigma\)^{2/k}\right\}$ with some
appropriate $\alpha>0$.
 
In the case $x\ge\sqrt{2C_k}n^{k/2}\sigma^{k+1}$ let us
choose $M=2^m$ in such a way that $\dfrac1{2^{(k+3)}}\le
\dfrac{C_kM^{k+1}}{nx^2}\le \dfrac14$. (There is a number $M=2^m$
satisfying this inequality, because by our temporary conditions
$nx^2\ge 2^{k+3}C_k$, i.e.\ $\dfrac{C_kM^{k+1}}{nx^2}\le
\dfrac1{2^{k+3}}$ for $M=1$.) With the above
choice of $M$ the inequality $\dfrac M{n\sigma^2}\ge
2^{-(k+2)/(k+1)}$ holds, since $\dfrac1{2^{k+3}}\le\dfrac{C_k
M^{k+1}}{nx^2}\le\dfrac 12\(\dfrac M{n\sigma^2}\)^{k+1}$ for
$x\ge\sqrt{2C_k}n^{k/2}\sigma^{k+1}$. Hence, by formula (3.20)
and the inequality $\dfrac1{2^{(k+3)}}\le
\dfrac{C_kM^{k+1}}{nx^2}\le \dfrac14$
$$
\align
P\(|J_{n,k}(f)|>x\)&\le
\(\frac{C_k\sigma^2
M^k}{x^2}\)^M\(\dfrac{2^{(k+2)/(k+1)}M}{n\sigma^2}\)^M
=\(\frac{2^{(k+2)/(k+1)}C_k M^{k+1}}{nx^2}\)^M\\
&\le\(\frac{2^{(k+2)/(k+1)}}{4}\)^M
\le \exp\left\{-\alpha(nx^2)^{1/(k+1)} \right\}
\endalign
$$
with some appropriate $\alpha>0$.
 
The above calculations show that formula (1.3) holds if
$\dfrac x\sigma\ge 2^{(k+1)/2}C_k^{1/2}$ and $nx^2\ge 2^{k+3}C_k$.
By choosing the constant $C>0$ sufficiently large in formula (1.3)
we can achieve that the right-hand side of (1.3) is greater than 1
if one of the above inequalities does not hold. With such a choice
of the constant $C$ we get that formula (1.3) holds for all $x\ge0$.
 
Now we turn to the proof of Theorem 2. \medskip\noindent
{\it The proof of Theorem 2.}\/ Let us first consider the case when
the measure $\mu$ is non-atomic. Let us take a canonical function
$f(x_1,\dots,x_k)$ of $k$ variables. In this case we can write
$$
\align
&n^{-k/2}I_{n,k}(f)=\frac{n^{k/2}}{k!}\int'
f(u_1,\dots,u_k)\mu_n(\,du_1)\dots\mu_n(du_k)\\
&\qquad=\frac{n^{k/2}}{k!}\int'
f(u_1,\dots,u_k)(\mu_n(\,du_1)-\mu(\,du_1))\dots(\mu_n(du_k)-\mu(du_k))
=J_{n,k}(f).       \tag3.21
\endalign
$$
(Let us remark that we have exploited at this point that the measure
$\mu$ is non-atomic. The canonical property of the function $f$ implies
that $\int f(u_1,\dots,u_k)\mu(\,du_j)=0$ for any $1\le j\le k$ and
$u_s\in X$, $j\neq s$. But we need the fact that this relation remains
valid if we omit the points $u=u_j$, $j\neq s$, from the domain of
integration.) Theorem 1 and formula (3.21) imply Theorem~2 in the case
when $\mu$ is non-atomic.
 
The general case can be reduced to the case of non-atomic measure
$\mu$ with the help of the following construction. Given a space
$(X,\Cal X)$ with a measure $\mu$ on it together with a sequence of
independent random variables $\xi_1,\dots,\xi_n$ and a canonical
function $f(x_1,\dots,x_k)$ on the product space $(X^k,\Cal X^k)$
consider the space $(Y,\Cal B,\lambda)$ where $Y$ is the unit
interval, $\Cal B$ is the Borel $\sigma$-algebra and $\lambda$ is the
Lebesgue measure on it. Then define the product space $(\bar
X,\bar{\Cal X},\bar\mu)$ as $\bar X=X\times Y$, $\bar{\Cal X}=\Cal
X\times \Cal A$ and $\bar\mu=\mu\times \lambda$. Define the function
$\bar f$ on the $k$-th power of this new space $(\bar X^k,\bar {\Cal
X}^k)$ as $\bar f((x_1,y_1),\dots,(x_k,y_k))=f(x_1,\dots,x_k)$, i.e.\
the function $\bar f$ does not depend on the new coordinates. Finally,
let us consider a sequence of independent random variables
$\eta_1,\dots,\eta_n$ with uniform distribution on the unit interval
$[0,1]$ which are independent also of the random variables
$\xi_1,\dots,\xi_n$. Define the random variables
$(\bar\xi_1,\dots,\bar\xi_n)=((\xi_1,\eta_1),\dots,(\xi_n,\eta_n))$,
and let $\bar\mu_n$ denote their empirical distribution function. Then
$I_{n,k}(f)=I_{n,k}(\bar f)$, and the functions $f$ and $\bar f$ are
simultaneously canonical. This means that the result of Theorem~2 can
be applied for the function $\bar f$, hence it also holds for the
function $f$.
 
\beginsection 4. On some basic properties of the random integrals
$J_{n,k}$
 
First we prove Lemma~1 which expresses the expected value of the
random integrals $J_{n,k}(f)$. \medskip\noindent
{\it The proof of Lemma~1.} Let us prove formula (3.8) first in the
special case when $f(x_1,\dots,x_k)=I_{A_1}(x_1)\cdots I_{A_k}(x_k)$,
where $A_1,\dots,A_k$ are disjoint measurable subsets of $X$ and
$I(A)$ denotes the indicator function of the set $A$. Let us introduce
the (random) signed measures $\nu_l$, $1\le l\le n$, on $(X,\Cal X)$
$$
\nu_l(A)=\nu_l(A,\oo)=\left\{
\aligned
1-\mu(A)&\qquad\text{if }\xi_l(\oo)\in A \\
-\mu(A)&\qquad\text{if }\xi_l(\oo)\notin A
\endaligned \right.
\qquad\text{for $A\in \Cal X$ and } 1\le l\le n.
$$
Then $\mu_n-\mu=\dfrac1n\summ_{l=1}^n\nu_l$ and
$$
EJ_{n,k}\(I_{A_1}(x_1)\cdots I_{A_k}(x_k)\)=\frac{n^{-k/2}}{k!}
E\(\prod_{m=1}^k \(\summ_{l=1}^n\nu_l(A_m)\)\). \tag4.1
$$
To calculate the expression in (4.1) first we fix some set of
indices $D=\{m_1,\dots,m_r\}\subset\{1,\dots,k\}$ of cardinality~$r$
together with a number $1\le l\le n$ and show that
$$
E\(\prod_{t=1}^r\nu_l(A_{m_t})\)=(r-1)(-1)^{r-1}
\prod_{t=1}^r\mu(A_{m_t}). \tag4.2
$$
Indeed, let us consider separately the cases when $\xi_l\in A_{m_t}$
for an $m_t\in D$ and the case when $\xi_l\notin A_{m_t}$
for all $m_t\in D$, and write
$$
\align
&P(\nu_l(A_{m_t})=1-\mu(A_{m_t}))\\
&\quad=P(\nu_l(A_{m_t})=1-\mu(A_{m_t}),\;\nu_l(A_{m_q})=-\mu(A_{m_q})
\text{ if } 1\le q\le r \text{ and } q\neq t)=\mu(A_{m_t})
\endalign
$$
for all $1\le t\le r$, and $P\(\nu_l(A_{m_t})=-\mu(A_{m_t}) \text
{ for all } 1\le t\le r)\) =1-\summ_{t=1}^r\mu(A_{m_t})$. Hence
$$
\align
E\(\prod_{t=1}^r\nu_l(A_{m_t})\)&=\sum_{t=1}^r \mu(A_{m_t})
\(\prod_{1\le q\le r,\;q\neq t}(-\mu(A_{m_q}))+
\prod_{q=1}^r(-\mu(A_{m_q}))\)\\
&\quad+\(1-\sum_{t=1}^r\mu(A_{m_t})\)
\prod_{q=1}^r(-\mu(A_{m_q}))
=(r-1)(-1)^{r-1}\prod_{t=1}^s\mu(A_{m_t}).
\endalign
$$
Take a partition $D_{1},\dots,D_{s}$ of the set $\{1,\dots,k\}$
such that $|D_{u}|=r_u$, $1\le u\le s$, with some prescribed numbers
$r_u$ and some integers $l_u$, $1\le l_u\le n$, $1\le u\le s$, such
that $l_u\neq l_{u'}$ if $u\neq u'$. The random vectors
$(\nu_{l_u}(A_1),\dots, \nu_{l_u}(A_k))$ are independent for
$u=1,\dots,s$. This fact together with relation (4.2) imply that
$$
\aligned
E\(\prodd_{u=1}^s\prodd_{v\in D_u}\nu_{l_u}(A_{v})\)
&=\prodd_{u=1}^s E\(\prodd_{v\in D_u}\nu_{l_u}(A_{v})\)
=\prod_{u=1}^s(r_u-1) (-1)^{r_u-1}\prod_{j=1}^k\mu(A_j)\\
&=(-1)^{k-s}\prod_{u=1}^s(r_u-1)
\int f(u_1,\dots,u_k)\mu(\,du_1)\dots\mu(\,du_k)
\endaligned \tag4.3
$$
with $f(u_1,\dots,u_k)=I_{A_1}(u_1)\cdots I_{A_k}(u_k)$. We get
formula (3.8) together with the constant defined in (3.9) in the
special case we have considered by means of relations (4.1) and
(4.3) if we observe that by carrying out the multiplication at the
right-hand side of (4.1) we get a sum with terms of the forms defined
in (4.3). (To get a term of the form given in (4.3) we prescribe
which measures $\nu_{u_l}$, $\le u_l\le n$, $1\le l\le k$, we take,
together with the partition of the set $\{1,\dots,k\}$ which tells
for which indices $l$ the measure $\nu_{u_l}$ is chosen. We can
choose $\binom ns s!B(r_1,\dots,r_s)$ such terms, where
$B(r_1,\dots,r_s)$ denotes the number of partitions of the set
$\{1,\dots,k\}$ to subsets of size $r_1$, $r_2$,\dots and $r_s$.)
 
Relation (3.8) also holds for any linear combination of functions of
the above type. Since the set of such linear combinations is dense in
the $L_1$ norm in the space $(X^k,\Cal X^k,\mu^k)$ to prove formula
(3.8) for general functions $f$ it is enough to show that
$$
E |J_{n,k}(f)|\le C(n,k) \|f\|_1\quad \text{for all functions $f$ of
$k$ variables} \tag4.4
$$
with an appropriate constant $C(k,n)<\infty$. Let us emphasize that
since in the problem considered here the dimension $k$ of the
function $f$ and the size $n$ of the sample $\xi_1,\dots,\xi_n$ are
fixed an arbitrary large constant $C(n,k)$ would suffice in (4.4).
 
Moreover, relation (4.4) can be simplified. By carrying out
the multiplications of the measures $\mu_n(\,du_j)-\mu(\,du_j)$
in formula (1.1) we get a sum with finitely many terms, and it is
enough to prove the bound given in (4.4) for each of these terms.
Moreover, by carrying out the integration with respect the
deterministic measures $\mu(\,du_j)$ we get that it is
enough to show that
$$
\align
E\int' &|f(u_1,\dots,u_{\bar k})|\mu_n(\,du_1)\dots\mu_n(\,du_{\bar k})
\le C(n,k) \|f\|_1 \\
&\qquad \text{for all functions $f$ of
$\bar k\le k$ variables}. \tag4.5
\endalign
$$
Even this relation can be reduced with the help of the triangular
inequality in the $L_1$ norm and the fact that the linear
combinations of the indicator functions of measurable sets are
dense in the space $L_1(X^{\bar k},\Cal X^{\bar k},\mu^{\bar k})$.
It is enough to prove the following simpler statement: For any $\bar
k\le k$ and measurable set $A\in (X^{\bar k},\Cal X^{\bar k})$
$$
E\int'_A \mu_n(\,du_1)\dots\mu_n(du_k)\le
C(n,k) \mu^{\bar k}(A). \tag4.6
$$
This statement is obvious, since we can express the integral at
the left-hand side of (4.6) as
$n^{\bar k}E\int'_A \mu_n(\,du_1)\dots\mu_n(du_k)=
n(n-1)\cdots(n-\bar k+1)\mu^{\bar k}(A)$, and this implies
relation~(4.6). The above identity holds, since to calculate the
integral at the left-hand side of (4.6) we have to count the number
of such $\bar k$ tuples $(\xi_{j_1},\dots,\xi_{j_{\bar k}})$ for
which $(\xi_{j_1},\dots,\xi_{j_{\bar k}})\in A$ and apply the
identity $P\((\xi_{j_1},\dots,\xi_{j_{\bar k}})\in A\)=\mu^{\bar
k}(A)$. (Let us remark that at this point we heavily exploited that
in the definition of the random integrals $J_{n,k}$ we deleted the
diagonals from the domain of integration.)
 
To complete the proof of Lemma 1 it is enough to check inequality
(3.10). For this aim let us observe that because of the factors
$r_u-1$ in (3.9) we may assume in this formula that $r_u\ge 2$ for
all $1\le u\le s$, hence $s\le\[\frac k2\]$, where $[x]$ denotes
the integer part of the number $x$. Beside this, we shall apply the
rather rude inequality $\prodd_{u=1}^s (r_u-1)\le 4^k$ if $r_u\ge1$,
$1\le u\le s$ and $\summ_{u=1}^s r_u=k$. Thus we get that
$$
\aligned
|B_{n,k}|&=\left|\frac1{k! n^{k/2}}\sum_{s=1}^k (-1)^{k-s}
\binom ns s! \!\!\!
\sum\Sb (r_1,\dots,r_s)\\ r_j\text{ is integer } r_j\ge1, \;1\le j\le
s\\ s_1+\cdots+r_s=k\endSb  \!\!\!
(r_1-1)\cdots(r_s-1)B(r_1,\dots,r_s) \right|\\
&\le \frac{4^k}{k! n^{k/2}}\sum_{s=1}^{\[\frac k2\]} n(n-1)\cdots(n-s+1)
\sum\Sb (r_1,\dots,r_s)\\ r_j\text{ is integer } r_j\ge1, \;1\le j\le
s\\ s_1+\cdots+r_s=k\endSb B(r_1,\dots,r_s).
\endaligned \tag4.7
$$
To estimate the above expression $B(n,k)$ let us first show that the
inner sum at the right-hand side of this expression satisfies the
inequality
$$
\sum\Sb (r_1,\dots,r_s)\\ r_j\text{ is integer } r_j\ge1, \;1\le j\le
s\\ s_1+\cdots+r_s=k\endSb B(r_1,\dots,r_s) \le
2^ks^{k-s},\quad\text{for all }1\le s\le k. \tag4.8
$$
Indeed, the expression at the left-hand side of (4.8) equals the
number of partitions of the set $\{1,\dots,k\}$ into $s$ non-empty
sets. This can be bounded with the help of the following construction:
Choose a set $\{u_1,\dots,u_s\}$, $1\le u_1<u_2<\cdots u_s\le k$,
what can be done in $\binom ks$ way. Let us construct a partition
$D_1,\dots,D_s$ of the set $\{1,\dots,k\}$ so that $u_j\in D_s$,
$1\le j\le s$, and distribute the remaining $k-s$ points of
$\{1,\dots,k\}$among the sets $D_j$, $1\le j\le s$. This can be done
in $s^{k-s}$ way. Beside this $\binom ks\le 2^k$. We have defined
no more than $2^ks^{k-s}$ partitions. I claim that we have constructed
all possible partitions of $s$ elements. Indeed, let us put the
elements of a partition $D_1,\dots,D_s$ in such an order that $j<j'$
if the smallest element $\bar u_j$ of $D_j$ is smaller than the
smallest element $\bar u_{j'}$ of $D_{j'}$. This partition
will be obtained if the set $\{\bar u_1,\dots,\bar u_s\}$
is chosen in the first step of the above construction, and then the
remaining terms are appropriately distributed among the sets of
partition.
 
Relations (4.7), (4.8) and the condition $1\le k\le\frac n2$ yield that
$$
|B_{n,k}|\le \frac{8^k}{k! n^{k/2}}\sum_{s=1}^{\[\frac k2\]} n^s
s^{k-s} \le \frac{\bar C^k}{k^k}\sum_{s=1}^{\[\frac k2\]}
n^{s-k/2} s^{k-s} \le \frac{(2\bar C)^k}{k^{k/2}}\sum_{s=1}^{\[\frac
k2\]}\(\frac sk\)^{k-s} \le\frac {C^k}{k^{k/2}}
$$
with some appropriate constants $\bar C>0$ and $C>0$. This means that
formula (3.10) holds. Lemma~1 is proved.
\medskip
Now we turn to the proof of the diagram formula described in Lemma~2.
Let me remark that similar results were already proved for
random integrals with respect to a Poisson process, (see~[8]).
These two problems and their possible proofs are similar, although
in the present problem some additional technical difficulties appear
because of the absence of independence for the empirical measure
$\mu_n$ of disjoint sets.
\medskip\noindent
{\it The proof of Lemma 2.}\/ Let us consider first the case when
$f(x_1,\dots,x_{k_1})=\prodd_{j=1}^{k_1} I_{A_j}(x_j)$,
$g(x_1,\dots,x_{k_2})=\prodd_{j'=1}^{k_2} I_{B_{j'}}(x_{j'})$,
the sets $A_1,\dots,A_{k_1}$ are disjoint and the same relation
holds for the sets $B_1,\dots,B_{k_2}$. (On the other hand, the sets
$A_j\cap B_{j'}$ may be non-empty. We denote again the indicator
function of a set $A$ by $I_A(x)$.) Let us choose a parameter $\e>0$
which we shall let go to zero and consider a partition
$A_j=\bigcupp_{s=1}^{u(j)}A(j,s)=\bigcupp_{s=1}^{u(j,\e)}A(j,s)$,
$1\le j\le k_1$, $B_{j'}=\bigcupp_{t=1}^{v(j')}B(j',t)=
\bigcupp_{t=1}^{v(j',\e)}B(j',t)$, $1\le j'\le k_2$, of the sets
$A_j$ and $B_{j'}$ such that $\mu(A(j,s))\le \e$,
$\mu(B(j',t))\le \e$, for all $1\le j\le k_1$, $1\le j'\le k_2$, $1\le
s\le u(j)$, $1\le t\le v(j')$, and beside this either $A(j,s)=B(j',t)$
or $A(j,s)\bigcap B(j',t)=\emptyset$ for all pairs $(j,s)$ and
$(j',t)$. We can write
$$
\align
J_{n,k_1}(f)J_{n,k_2}(g)=\frac{n^{(k_1+k_2)/2}}{k_1!k_2!} \!\!
\sum\Sb (s_1,\dots,s_{k_1})\:\\ 1\le s_j\le u(j),\\
1\le j\le k_1\endSb \!
&\sum\Sb (t_1,\dots,t_{k_2})\:\\ 1\le t_{j'}\le v(j'),\\
1\le j'\le k_2\endSb  \!
\prod_{j=1}^{k_1}\(\mu_n(A(j,s_j))- \mu(A(j,s_j))\)\\
&\prod_{j'=1}^{k_2}\(\mu_n(B(j',t_{j'}))-
\mu(B(j',t_{j'}))\). \tag4.9
\endalign
$$
Let $\Cal B(\Cal N)=\Cal B(\Cal N,k_1,k_2)$ denote the set of all
diagrams with first row $\{1,\dots,k_1\}$ and second row
$\{1,\dots,k_2\}$. We rewrite formula (4.9) with the help of
these diagrams by putting the products in the sum in (4.9) into certain
groups in dependence on the behaviour of the indices defining them.
Put
$$
J_{n,k_1}(f)J_{n,k_2}(g)=\frac{n^{(k_1+k_2)/2}}{k_1!k_2!}
\sum_{B(\Cal N)\in \Cal B(\Cal N)} Z_{\Cal N} \tag4.10
$$
with
$$
\aligned
Z_{\Cal N}=Z_{\Cal N}(\e)=
\sumn\Sb (s_1,\dots,s_{k_1}),\,
(t_1,\dots,t_{k_2})\:\\
1\le s_j\le u(j),\,1\le j\le k_1,\\
1\le t_{j'}\le v(j'),\,1\le j'\le k_2, \endSb
&\prod_{j=1}^{k_1}\(\mu_n(A(j,s_j))- \mu(A(j,s_j))\)\\
&\qquad\prod_{j'=1}^{k_2}\(\mu_n(B(j',t_{j'}))-
\mu(B(j',t_{j'}))\),
\endaligned \tag4.11
$$
where $\sumn$ means that summation is taken for such pairs of vectors
$(s_1,\dots,s_{k_1})$ and $(t_1,\dots,t_{k_2})$ for which
$A(j,s_j)=B(j',t_{j'})$ if $(j,j')$ is an edge of the diagram
$B(\Cal N)$ and $A(j,s_j)\cap B(j',t_{j'})=\emptyset$ if $(j,j')$
is not an edge of the diagram $B(\Cal N)$. If there are no vectors
$(s_1,\dots,s_{k_1})$ and $(t_1,\dots,t_{k_2})$ with such a property,
then we define $Z_{\Cal N}=0$.
 
Let us consider a diagram $B(\Cal N)\in \Cal B(\Cal N)$ with edges
$\Cal N=\{(j_1,j_1'),\dots,(j_l,j_l')\}$ and consider an
approximation $\bar Z_{\Cal N}$ of $Z_{\Cal N}$ which we obtain in
the following way: If $(j,j')\in \Cal N$ then we replace the
term $(\mu_n(A(j,s_j))-\mu(A(j,s_j)))
(\mu_n(B(j',t_{j'}))-\mu(B(j',t_{j'})))
=\(\mu_n(A(j,s_j))-\mu(A(j,s_j))\)^2$ by $\dfrac1n\mu_n(A(j,s))$ and
do not change the other terms. (The following heuristic argument may
lead to the introduction of such an approximation. The expression
$\(\mu_n(A(j,s_j))-\mu(A(j,s_j))\)^2$ can be well approximated by
$\mu_n(A(j,s_j))^2$, the other terms being negligible. On the other
hand, $\mu_n(A(j,s_j))^2\sim\frac1n \mu_n(A(j,s_j))$, since the
probability of the event that at least two sample points $\xi_j$
fall in a set $A(j,s_j)$ of small $\mu$ measure is negligibly small.
Hence $\mu_n^2(A(j,s_j))-\frac1n \mu_n(A(j,s_j))=0$ with probability
almost~1.) More explicitly, we define $\bar Z_{\Cal N}$ for all
$B(\Cal N)\in\Cal B(\Cal N)$ as
$$
\align
\bar Z_{\Cal N}=\bar Z_{\Cal N}(\e)=
\sumn\Sb (s_1,\dots,s_{k_1}),\,
(t_1,\dots,t_{k_2})\:\\
1\le s_j\le u(j),\,1\le j\le k_1,\\
1\le t_{j'}\le v(j'),\,1\le j'\le k_2\endSb
&\prod_{j\: j\in
\{1,\dots,k_1\}\setminus\{j_1,\dots,j_l\}}\(\mu_n(A(j,s_j))-
\mu(A(j,s_j))\)\\
&\qquad\qquad\qquad \prod_{m=1}^l\(\frac
1n\mu_n(A(j_m,s_{j_m}))\)\\
&\hskip-4truecm\prod_{j'\: j'\in \{1,\dots,k_2\}
\setminus\{j_1',\dots,j_l'\}}
\(\mu_n(B(j',t_{j'}))-\mu(B(j',t_{j'}))\), \tag4.12
\endalign
$$
 
I claim that for all $B(\Cal N)\in\Cal B(\Cal N)$ the inequality
$$
\left| Z_{\Cal N}(\e)-\bar Z_{\Cal N}(\e)\right|\le C(k_1+k_2,n)
\e \tag4.13
$$
holds with some constant $C(k_1+k_2,n)<\infty$ which does not depend
on $\e$. We do not need a good estimate for this constant. What is
important for us is that relations (4.10) and (4.13) have the
consequence
$$
J_{n,k_1}(f)J_{n,k_2}(g)=\lim_{\e\to0}
\frac{n^{(k_1+k_2)/2}}{k_1!k_2!}
\sum_{B(\Cal N)\in \Cal B(\Cal N)}\bar Z_{\Cal N}(\e). \tag4.14
$$
Before proving relation (4.13) let us show that formulas (4.14) and
(4.12) imply relations (3.11) and (3.12) in the special case we are
considering. To do this we express $\bar Z_{\Cal N}(\e)$
(which as we shall see actually does not depend on the parameter
$\e$) in an appropriate form. To do this let us write the terms in the
middle product of (4.12) in the form $\frac 1n
\mu_n(A(j_m,s_{j_m}))=\frac1n \((\mu_n(A(j_m,s_{j_m}))
-\mu(A(j_m,s_{j_m}))) +\mu(A(j_m,s_{j_m}))\)$, and let us rewrite the
terms $\bar Z_{\Cal N}$ as
$$
\align
\bar Z_{\Cal N}&=n^{-l}\int
\prod_{t=1}^l R_{j_t,j'_t} f\circ g(x_1,\dots,x_{k_1+k_2})
\prod_{j\: j\in \{1,\dots,k_1\}\setminus\{j_1,\dots,j_l\}}
\(\mu_n(d x_j)-\mu(\,dx_j)\)   \\
&\qquad \prod _{j'\: j'\in \{1,\dots,k_2\}
\setminus\{j_1',\dots,j_l'\}} \(\mu_n(d x_{j'})-\mu(\,dx_{j'})\)
\tag4.15   \\
&\qquad\qquad \prod_{s=1}^l \((\,d\mu_n(x_{j_s})-\,d\mu(x_{j_s}))
+\,d\mu(x_{j_s})\),
\endalign
$$
where $(j_1,j'_1)$,\dots $(j_l,j'_l)$ are the edges of the diagram
$B(\Cal N)$, and the function $f\circ g$ and operators $R_{j_t,j'_t}$
are defined in formulas (3.4) and (3.3).

Given a diagram $B(\Cal N)$ let us introduce the class of coloured
diagrams $\Cal C(B(\Cal N))$ consisting of all coloured diagrams
$B(\Cal N,\Cal N_1)$ we obtain by colouring the edges of the diagram
of $B(\Cal N)$ by $+1$ or $-1$. By carrying out the multiplications
in formula (4.15) by expressing the last line of this expression as
the sum of $d\mu_n(x_{j_s})-\,d\mu(x_{j_s})$ and $\,d\mu(x_{j_s})$
we can rewrite $\bar Z_{\Cal N}$ with the help of the coloured diagrams
in $\Cal C(B(\Cal N))$ as
$$
\bar Z_{\Cal N}=\sum_{p=1}^l\sum\Sb B(\Cal N,\Cal N_1)\in\Cal C(B(\Cal
N)),\\  \text{and }  |\Cal N_1|=p \endSb
n^{-l} \dfrac{(k_1+k_2-l-p)!}{n^{(k_1+k_2-l-p)/2}} J_{n,k_1+k_2-l-p}
(\overline {f\circ g}_{B(\Cal N,\Cal N_1)})
$$
The last relation together with formula (4.14) yield that
$$
J_{n,k_1}(f)J_{n,k_2}(g)=\!\!\!\! \sum_{l=0}^{\min (k_1,k_2)} \!\!
\sum_{p=0}^l \sum_{B(\Cal N,\Cal N_1)\in B(l,p)} \!\!\!\!
 \dfrac{(k_1+k_2-l-p)!}{n^{(l-p)/2}k_1!k_2!} J_{n,k_1+k_2-l-p}
(f\circ g_{B(\Cal N,\Cal N_1)}).
$$
This relation is equivalent to the diagram formula described in
formulas (3.11) and (3.12). Hence to prove Lemma 2 in the special
case considered now it is enough to prove formula (4.13).
 
First we formulate an inequality and show that formula (4.13)
follows from it. To formulate it let us introduce the following
five (random) functions $\rho^{(m)}$, $1\le m\le 5$, defined on the
measurable sets $A\in \Cal X$ of the space $(X,\Cal X)$. Put
$\rho^{(1)}(A)=\mu_n(A)$, $\rho^{(2)}(A)=-\mu(A)$,
$\rho^{(3)}(A)=-2\mu(A)\mu_n(A)$, $\rho^{(4)}(A)=\mu(A)^2$ and
$\rho^{(5)}(A)=\mu_n(A)^2-\frac1n\mu_n(A)$ for all $A\in \Cal X$.
Let us fix some positive integer $r$ and $r$ disjoint subsets
$C_j\in\Cal X$ of the space $X$, $1\le j\le r$, such that
$\mu(C_j)\le\e$ for all $1\le j\le r$  with some $r\ge1$. Let us
choose some sequence of integers $m(j)$, $1\le j\le r$, such that
$1\le m(j)\le 5$ for all $1\le j\le r$ and $m(j)\ge3$ for at least
one of these indices. I claim that
$$
E\prod_{j=1}^r\left|\rho^{m(j)}(C_j)\right|\le  \e\bar C(r,n)
\prod_{j=1}^r\mu(C_j) \tag4.16
$$
for the above defined functions $\rho^{(m)}$ with some appropriate
constant $\bar C(r,n)$. (One could get rid of the dependence of
$n$ in the constant $\bar C(r,n)$ but such an improvement of the
inequality (4.16) has no great importance for us.)
 
First I show that relation (4.13) can be proved with the help
of the estimate (4.16). For this aim I rewrite the expressions
$Z_{\Cal N}$ and $\bar Z_{\Cal N}$ with the help of the above defined
functions $\rho^{(m)}$. To get a good representation let us first
consider those pairs $(j,j')$ for which $(j,j')\in \Cal N$ and
write the terms of the form $(\mu_n(A(j,s_j))-\mu(A(j,s_j)))
(\mu_n(B(j',t_{j'}))-\mu(B(j',t_{j'})))=\(\mu_n(A(j,s_j))
-\mu(A(j,s_j))\)^2$ in the expression (4.11) with such indices $j$
and $j'$ by means of these new functions as
$\rho^{(5)}(A(j,s_j))+\frac1n\rho^{(1)}(A(j,s_j))
+\rho^{(3)}(A(j,s_j))+\rho^{(4)}(A(j,s_j))$. Similarly, the terms
$\frac1n\mu_n(A(j,s))$ in the expression (4.12) can be rewritten as
$\frac1n\rho^{(1)}(A(j,s_j))$. If $j$ is not an end-point of an edge in
the first row of the diagram $\Cal N$, then let us rewrite the terms of
the form $\mu_n(A(j,s_j))-\mu(A(j,s_j))$ as
$\rho^{(1)}(A(j,s_j))+\rho^{(2)}(A(j,s_j))$ both in formula (4.11)
and (4.12). Similarly, if $j'$ is not an end-point of an edge in the
 second row of the diagram $\Cal N$, then let us rewrite the terms of
the form $\mu_n(B(j',t_{j'}))-\mu(B(j',t_{j'}))$ as
$\rho^{(1)}(B(j',t_{j'}))+\rho^{(2)}(B(j',t_{j'}))$ in these formulas.
 
Let us consider the difference $Z_{\Cal N}-\bar Z_{\Cal N}$
rewritten in the above form as the function of the terms
$\rho^{m(j)}(A_j,s_j)$, $1\le j\le k_1$ and
$\rho^{m(j')}(B_{j'},t_{j'})$, $j'\in \{1,\dots,k_2\}
\setminus \{j_1',\dots,j'_l\}$. In such a way we can express $Z_{\Cal
N}-\bar Z_{\Cal N}$ as the linear combination of some terms with
coefficients between zero and 1 which have the form
$$
\prodd_{s=1}^{k_1}\rho^{m(j)}(A(j,s_j))\prodd_{j'\in \{1,\dots,k_2\}
\setminus\{j_1',\dots,j'_l\}}\rho^{m(j')}(B(j',t_{j'})), \tag4.17
$$
where $\{j_1',\dots,j'_l\}$ is the set consisting of the end-points
of the edges of the diagram $B(\Cal N)$ in the second row.
Beside this, it follows from the above construction that in the terms
described in (4.17) $m(j)\ge3$ for at least one of the
indices $j_1,\dots,j_l$, where $j_1,\dots,j_l$ are the end-points of
the edges of the diagram $\Cal N$ in the first row. Since by our
assumption the inequalities $\mu(A(j,s))\le \e$ and
$\mu(B_{j'},t)\le\e$ hold for all
terms we consider, formula (4.16) can be applied with the choice
$r=k_1+k_2-l$ to estimate all terms of the form (4.17). It yields that
$$
\aligned
&E\left|\prodd_{s=1}^{k_1}\rho^{m(j)}(A(j,s_j))\prodd_{j'\in
\{1,\dots,k_2\}
\setminus\{j_1',\dots,j'_l\}}\rho^{m(j')}(B(j',t_{j'}))\right|\\
&\qquad \le \e \bar C(k_2,n) \prodd_{s=1}^{k_1}\mu(A(j,s_j))
\prodd_{j'\in \{1,\dots,k_2\}\setminus
\{j_1',\dots,j'_l\}}\mu(B(j',t_{j'})) \\
&\qquad=\e\bar C(k_1+k_2,n) \mu^{k_1+k_2-l}\(\prodd_{s=1}^{k_1}
A(j,s_j)\times  \prodd_{j'\in \{1,\dots,k_2\}
\setminus\{j_1',\dots,j'_l\}}   B(j',t_{j'})\),
\endaligned \tag4.18
$$
where $\mu^{k_1+k_2-l}$ denotes the $k_1+k_2-l$-fold direct product
of the measure $\mu$ with itself on the space $(X^{k_1+k_2-l},\Cal
X^{k_1+k_2-l})$, and the product at the right-hand side of (4.18)
denotes the direct product of the corresponding sets.
 
The estimate (4.13) follows from inequality (4.18) if we sum up
this estimate for all such terms of the form (4.17) which appear in
the difference $Z_{\Cal N}-\bar Z_{\Cal N}$ when it is written in the
above described form. We can see this if we observe that any two sets
$A(j,s)$ and $B(j',t)$ appearing in one of these expressions are
disjoint unless they are both functions of the same type
$A(\cdot,\cdot)$ or $B(\cdot,\cdot)$, and their arguments $(j,s_j)$
or $(j',t_{j'})$ agree. This fact implies that those sets
$\prodd_{s=1}^{k_1}A(j,s_j)\times \prodd_{j'\in
\{1,\dots,k_2\} \setminus\{j_1',\dots,j'_l\}} B(j',t_{j'})$
at the right-hand side of (4.18) whose $\mu^{k_1+k_2-l}$ measures
have to be summed up to get an estimate
for the difference $E|Z_{\Cal N}-\bar Z_{\Cal N}|$ are either
disjoint or they coincide. Since each term appears less than
$5^{k_1+k_2}$ times in one of these estimates this means that
formula (4.16) implies formula (4.13) with
the choice $C(k_1+k_2,n)=5^{k_1+k_2}\bar C(k_1+k_2,n)$.
 
To prove formula (4.16) let us express the functions
$\rho^{(1)}(C_j)$, $\rho^{(3)}(C_j)$ and $\rho^{(5)}(C_j)$ with
the help of the measures $\bar\nu_l$ defined as
$$
\bar\nu_l(A)=\bar\nu_l(A,\oo)=\left\{
\aligned
1&\qquad\text{if }\xi_l(\oo)\in A \\
0&\qquad\text{if }\xi_l(\oo)\notin A
\endaligned \right.
\qquad\text{for all $A\in \Cal X$ and } 1\le l\le n.
$$
We can write $\rho^{(1)}(C_j)=\dfrac1n\summ_{l_j=1}^n
\bar\nu_{l_j}(C_j)$, $\rho^{(3)}(C_j)=-\dfrac2n\mu(C_j)
\summ_{l_j=1}^n\bar\nu_{l_j}(C_j)$, and
$$
\rho^{(5)}(C_j)=\frac1{n^2}\summ_{l_j=1}^n\bar\nu_{l_j}(C_j)-
\(\frac1n\summ_{l_j=1}^n\bar\nu_{l_j}(C_j)\)^2=-\frac1{n^2}
\summ_{1\le l_j,l_j'\le n,\;l_j\neq l_{j'}}
\bar\nu_{l_j}(C_j)\bar\nu_{l_{j'}}(C_j).
$$
Let us observe that only pairs $(l_j,l_{j'})$ with $l_j\neq l_{j'}$
are considered at the right-hand side of the last formula, since
$\bar \nu_l(C)-\bar\nu_l^2(C)=0$ for all $C\in\Cal X$.
 
Let us insert these formulas to the left-hand side of formula (4.16),
carry out the multiplications, and let us bound the absolute value
of the sum obtained in such a way as the sum of the absolute value
of these terms. All these terms are products of expressions of the
form $\bar\nu_{l_j}(C_j)$, $\bar\nu_{l_{j'}}(C_j)$ or $\mu(C_j)$.
Let us observe that only those products have to be considered for
which the indices $l_j$ and $l_{j'}$ of the measures $\nu_{l_j}$
and $\nu_{l_{j'}}$ are all different for different indices $j$.
Otherwise the product we consider equals zero. We can see this by
observing that the sets $C_j$ are disjoint. Therefore if
$l_{j_1}=l_{j_2}$ for some $j_1\neq j_2$, then either
$\bar\nu_{l_{j_1}}(C_{j_1})=0$ or $\bar\nu_{l_{j_2}}(C_{j_2})=0$.
On the other hand, the measures $\bar\nu_l$ are independent for
different indices $l$, and this implies that the terms we have to
bound are the products of  $r$ independent random variables (some of
these variables may be constant). Since the sum we have to bound
contains less than $n^r$ terms to prove (4.16) it is enough to show
that the expectation of the absolute value of these terms can be
bounded by $\const \e\prodd_{j=1}^r\mu(C_j)$. Since $E\bar\nu_{l_j}(C_j)
=E\bar\nu_{l_{j'}}(C_j)=\mu(C_j)$, the product  whose expected
value has to be estimated is the product of $r$ such independent
random variables whose expected values can be bounded $\mu(C_j)$,
$1\le j\le r$. Moreover, there is at least one index $j$ for which
$m(j)\ge3$ in the expression (4.16). This property together with the
condition $\mu(C_j)\le\e$ imply that the product whose expected value
we have to estimate contains at least one factor whose expected
value is less than $2\e\mu(C_j)$. This means that the inequality needed
for the proof of relation (4.16) holds.
 
Let us observe that both sides of formula (3.11) is a bilinear form
of the functions $f$ and $g$. This fact implies that it holds not
only for those pairs $(f,g)$ of functions we considered before, but
also for any finite linear combinations of such functions $f$ and $g$.
Such linear combinations are dense in the spaces
$L_2(X^{k_1},\Cal X^{k_1},\mu^{k_1})$ and $L_2(X^{k_2},\Cal X^{k_2},
\mu^{k_2})$ respectively. Hence to complete the proof of Lemma~2 it is
enough to show that if a sequence of functions $f_n$ and $g_n$
satisfy the relations $\|f_n-f\|_2\to0$ and $\|g_n-g\|_2\to0$ with
some functions $f$ and $g$, then
$$
\lim_{n\to\infty}
E|J_{n,k_1}(f_n)J_{n,k_2}(g_n)-J_{n,k_1}(f)J_{n,k_2}(g)|=0 \tag4.19
$$
and
$$
\lim_{n\to\infty} E|J_{n,k_1+k_2-l-p}(f_n\circ g_n)_{B(\Cal N,\Cal
N_1)}-J_{n,k_1+k_2-l-p}(f\circ g)_{B(\Cal N,\Cal N_1)}|=0 \tag4.20
$$
for all coloured diagrams $B(\Cal N,\Cal N_1)\in B(l,p)$, $0\le p\le
l\le \min(k_1,k_2)$.
 
To prove relation (4.19) it is enough to prove that $E(J_{n,k_1}(f_n)
-J_{n,k_1}(f))^2=E((J_{n,k_1}(f_n-f))^2\to0$ if $\|f_n-f\|_2\to0$,
and $E(J_{n,k_2}(g_n-g))^2\to0$ if $\|g_n-g\|_2\to0$. Let us prove the
inequality $E(J_{n,k}(h))^2\le C(n,k_1)\|h\|_2^2$ with some constant
$C(k,n)$ for a function $h$ of $k$ variables. The relation we want to
prove is a consequence of this inequality. On the other hand, this
inequality can be deduced from some arguments presented in the proof
of Lemma~1. Indeed, similarly to the argument presented there we can
reduce the statement to be proved to the inequality $\int
\bar h^2(x_1,\dots,x_k)\mu_n(\,du_1)\dots\mu_n(\,du_{\bar k})
\le C(n,k)\|\bar h\|_2^2$ with an appropriate constant $C(n,k)$ for
any function $\bar h$ of $\bar k\le k$ variables such that
$\|\bar h\|_2^2<\infty$. But this statement is contained in
formula (4.5). We only have to apply it for the integrable function
$\bar h^2$.
 
Relation (4.20) can be deduced from the results proved at the end of
Lemma~1 in an even simpler way. Because of the results proved there it
is enough to check that $\left\|(f_n\circ g_n)_{B(\Cal N,\Cal N_1)}-
(f\circ g)_{B(\Cal N,\Cal N_1)}\right\|_1\to0$ if $\|f_n-f\|_2^2\to0$
and $\|g_n-g\|_2^2\to0$. On the other hand, it follows from the
Schwarz inequality and the contraction property of the
operator $P$ defined in (3.2) in the $L_2$ norm that
$$
\left\|(f\circ g_n)_{B\Cal N,\Cal N_1)}-(f\circ
g)_{B(\Cal N,\Cal N_1)}\right\|_1=
\left\|(f\circ (g_n-g))_{B(\Cal N,\Cal N_1)}\right\|_1\le\|f\|_2^2
\|g-g_n\|_2^2,
$$
hence $\left\|(f\circ g_n)_{B(\Cal N,\Cal N_1)}-(f\circ g)_{B(\Cal
N,\Cal N_1)}\right\|_1\to0$. Similarly,
$$
\left\|(f_n\circ g_n)_{B(\Cal N,\Cal N_1)}-(f_n\circ
g_n)_{B(\Cal N,\Cal N_1)}\right\|_1\to0.
$$
These inequalities imply relation (4.20). Lemma~2 is proved.
\medskip
Lemma 3 which gives an estimate on the second moment of a random
integral $J_{n,k}(f)$ can be proved as a consequence of Lemmas~1
and~2.
\medskip\noindent {\it The proof of Lemma 3.}\/ By Lemma 2 we have
$$
EJ_{n,k}(f)^2=\sum_{l=0}^k
\sum_{p=0}^l \frac{(2k-l-p)!}{(k-l)!^2(l-p)!p!}
n^{-(l-p)/2} E J_{n,2k-l-p}(f\circ f_{l,p})
$$
with the functions $f\circ f_{l,p}$ defined in formula (3.12). Hence
by formulas (3.8) and (3.10) in Lemma~1
$$
EJ_{n,k_1}(f)^2\le  \sum_{l=0}^k
\sum_{p=0}^l \frac{(2k-l-p)!}{(k-l)!^2(l-p)!p!}
n^{-(l-p)/2}\frac{C^{2k-l-p}}{(2k-l-p)^{(2k-l-p)/2}} \|f\circ
f_{l,p}\|_1 \tag4.21
$$
with some constant $C>0$. We claim that $\|f\circ f_{l,p}\|_1\le
\|f\|_2^2$. Since $f\circ f_{l,p}$ is the average of certain functions
of the form $f\circ f_{B(\Cal N,\Cal N_1)}$ with some diagram $B(\Cal
N,\Cal N_1)$ it is enough to show that
$\|f\circ f_{B(\Cal N,\Cal N_1)}\|_1\le \|f\|_2^2$ for an arbitrary
diagram $B(\Cal N,\Cal N_1)$. We get similarly to the proof of
formula (3.15) by considering first functions of the form
$\overline{f\circ f}_{B(\Cal N)}$ defined in formula (3.5) and
applying the Schwarz inequality for them that
$$
\left\|\overline {f\circ f}_{B(\Cal N)}\right\|_1 \le \int
f^2(x_1,\dots,x_k)\mu(\,dx_1) \dots \mu(\,dx_{k_1})=\|f\|_2^2
$$
for an arbitrary diagram $\Cal N$. Since the operators $P_j$ defined
in formula (3.2) are contractions also in $L_1$ norm, the last
inequality implies that the relation $\|f\circ f_{B(\Cal N,\Cal
N_1)}\|_1= \|\overline{f\circ f}_{B(\Cal N,\Cal N_1)}\|_1\le
\|\overline{f\circ f}_{B(\Cal N)}\|_1\le \|f\|_2^2$ holds, as we
claimed. By formula (4.21), the inequality $\|f\circ f_{l,p}\|_1\le
\|f\|_2^2$ and the Stirling formula we get that
$$
EJ_{n,k_1}(f)^2\le  C^k \sum_{l=0}^k
\sum_{p=0}^l \frac{(2k-l-p)^{(2k-l-p)/2}}{(k-l)^{2(k-l)}(l-p)^{l-p}p^p}
n^{-(l-p)/2} \|f\|_2^2 \tag4.22
$$
with a possibly new constant $C>0$.
 
Since the right-hand side of (4.22) contains less than $k^2$ terms
and $n\ge k$ this relation enables us to reduce the proof of Lemma~3
to the estimate
$$
(2k-l-p)^{(2k-l-p)/2} \le \bar C^k
(k-l)^{2(k-l)}(l-p)^{l-p}p^pk^{(l-p)/2} k^{-k}
$$
with an appropriate $\bar C\ge0$ for all $0\le p\le l\le k$. Since
$(k-l)^{k-l} p^p (l-p)^{l-p}\ge 3^{-k} k^k$ (this follows from
instance from the convexity of the function $x\log x$ if we take
logarithm), to prove the last inequality it is enough to show that
$$
\(\frac{2k-l-p}2\)^{(2k-l-p)/2} \le \tilde C^k
(k-l)^{(k-l)}k^{(l-p)/2} \quad\text{if } 0\le p\le l\le k
$$
with some $\tilde C>0$. By dividing both sides of the last inequality
by $k^{((2k-l-p)/2}$ and after this taking the $k$-th root of both
sides we get that the last inequality is equivalent to the statement
$$
\(\frac{2k-l-p}{2k}\)^{(2k-l-p)/2k} \le \tilde C
\(\frac {k-l}k\)^{(k-l)/k} \quad\text{if } 0\le p\le l\le k. \tag4.23
$$
Inequality (4.23) holds. This inequality holds since the left-hand side
can be bounded by $\supp_{0\le u\le1}u^u\le1$ from above, and the
right-hand side by $\tilde C\inff_{0\le u\le 1}u^u\ge \tilde Ce^{-1/e}$
from below. Lemma~3 is proved.
\medskip\noindent
{\it The proof of Lemma 4.}\/ Let us first investigate the main
result of Lemma 4 which states that the function $f\circ g$ is
$(r_1+r_2-(l-p),\sigma^2)$ dominated under appropriate conditions.
Let us fix a coloured diagram $B(\Cal N,\Cal N_1)=B(\Cal N,\Cal N_1,
k_1,k_2)$ with
$$
\Cal N=\{(j_1,j_1'),\dots,(j_l,j_l')\} \text{ and }\Cal N_1
=\{(j_{u(1)},j_{u(1)}'),\dots,(j_{u(p)},j_{u(p)}')\}, \tag4.24
$$
where $\{u(1),\dots,u(p)\}\subset\{1,\dots,l\}$. It is enough
to prove that the function $\overline{f\circ g}_{B(\Cal N,\Cal N_1)}$
defined in formula (3.6) is $(r_1+r_2-(l-p),\sigma^2)$ dominated.
We have certain freedom to change the order of the operators
$R_{j,j'}$ and $P_j$ in the definition of this function. By
exploiting this freedom we can rewrite the function $\overline
{f\circ g}_{B(\Cal N,\Cal N_1)}$ as
$$
\aligned
&\overline {f\circ g}_{B(\Cal N,\Cal N_1)}
(x_{s(1)},\dots,x_{s(k_1+k_2-l-p)}) \\
&\qquad=\prod_{s=1}^{l-p} R_{j_{v(s)},j'_{v(s)}}
\prod_{t=1}^p P_{j_{u(t)}} R_{j_{u(t)},j'_{u(t)}}
f\circ g(x_1,\dots,x_{k_1+k_2})
\endaligned \tag4.25
$$
with $\{v(1),\dots,v(l-p)\}=\{1,\dots,l\}\setminus
\{u(1),\dots,u(p)\}$, i.e.\ the operators $R_{j_{v(s)},j'_{v(s)}}$
are indexed by the edges
$$
\Cal M=\Cal N\setminus \Cal N_1=
\{(j_{v(1)},j_{v(1)}'),\dots,(j_{v(l-p)},j_{v(l-p)}')\} \tag4.26
$$
of the diagram $B(\Cal N,\Cal N_1)$. Beside this, we can choose the
order of the operators $R_{j_{v(s)},j'_{v(s)}}$ in an arbitrary way.
Let me recall that in our notations the operators written after
each other are applied in right to left order.
 
Lemma~4  will be proved with the help of an inductive procedure
described below. To carry it out let us introduce some notations.
Let us consider a function $F(x_{s(1)},\dots,x_{s(m)})$ with the
set of indices $S=\{s(1),\dots,s(m)\}\subset \{1,\dots,k_1+k_2\}$
together with a diagram $B(\Cal M)=B(\Cal M,S)$ consisting of two
rows $\{s(1),\dots,s(r)\}$ and $\{s(r+1),\dots, s(m)\}$, where the
index $r$ is defined by the relations $s(r)\le k_1<s(r+1)$. Beside
this, the diagram $B(\Cal M,S)$ has some edges $\Cal M=
\{(j_1,j'_1),\dots,(j_q,j'_q)\}$ with $0\le q\le \min(r,m-r)$,
$1\le j_1<\cdots j_q\le s(r)$,
$\{j_1,\dots,j_q\}\subset\{s(1),\dots,s(r)\}$,
$\{j'_1,\dots,j'_q\}\subset\{s(r+1),\dots, s(m)\}$ and
$j'_u\neq j'_{u'}$ if $u\neq u'$, $1\le u,u'\le q$. We
introduced the diagram $B(\Cal M,S)$ to indicate for which pairs
$(j,j')$ we want to apply the operator $R_{j,j'}$ for the function
$F$ during the inductive procedure leading to the proof of Lemma~4.
 
Let us call a partition $B_1,\dots, B_r$,
$B_u=\{s(d_{(u,1)}),\dots,s(d_{(u,v(u))})\}$, $1\le u\le r$, of the
set $S=\{s(1),\dots,s(m)\}$ together with a set of functions
$h_1,\dots,h_r$, such that the function $h_u=h_u\(x_{s(d_{(u,1)})},
\dots,x_{s(d_{(u,v(u))})}\)$ depends only on the variables indexed
by the elements of the set $B_u$ for all $1\le u\le r$, and
$\|h_u\|_\infty\le 1$, $\|h_u\|_2^2\le \sigma^2$, $1\le u\le r$, an
($r,\sigma^2)$ dominating system of the function $F$ if also the
relation
$$
|F(x_{s(1)},\dots,x_{s(m)})|\le \prod_{u=1}^r
h_u\(x_{s(d_{(u,1)})},\dots,x_{s(d_{(u,v(u))})}\)
\quad\text {for all } x_{s(l)}\in X,\; 1\le l\le m
$$
holds. Let us call this $(r,\sigma^2)$ dominating system regular
(with respect to the function $F$ and the diagram $B(\Cal M)$
attached to it) if the two end-points $j_i$ and $j'_i$ of the edges
of the diagram $B(\Cal M)$ are contained in different sets $B_{u(i)}$
and $B_{u'(i)}$ of the partition $B_1,\dots, B_r$ for all
$1\le i\le q$. We call an $(r,\sigma^2)$ dominating system
super-regular if all elements $B_1,\dots,B_r$ of the partition of
this system are contained either in the set $\{1,\dots,k_1\}$ or in
the set $\{k_1+1,\dots,k_1+k_2\}$. A super-regular $(r,\sigma^2)$
dominating system is clearly regular with respect to any diagram
$B(M)$, since all edges of a diagram $B(\Cal M)$ are going between
the vertices of the sets $\{1,\dots,k_1\}$ and
$\{k_1+1,\dots,k_1+k_2\}$.
 
It follows from the conditions of Lemma~4 that the function
$f\circ g(x_1,\dots,x_{k_1+k_2})$ defined in formula (3.4)
has a super-regular $(r_1+r_2,\sigma^2)$ dominating system. Indeed,
the function $f(x_1,\dots,x_{k_1})$ can be dominated by a partition
$B_1,\dots,B_{r_1}$ of the set $\{1,\dots,k_1\}$ together with some
functions $h_1,\dots,h_{r_1}$ such that the indices of the arguments
of the function $h_u$ are contained in $B_u$, $1\le u\le r_1$, and
the function $g\(x_{k_1+1},\dots,x_{k_1+k_2}\)$ by a partition
$B_{r_1+1},\dots,B_{r_1+r_2}$ of the set
$\{k_1+1,\dots,k_1+k_2\}$ together with some functions
$h_{r_1+1},\dots,h_{r_1+r_2}$ such that the indices of the arguments
of the function $h_u$ are contained in $B_u$, $r_1+1\le u\le r_1+r_2$.
Beside this $\|h_u\|_\infty\le1$ and $\|h_u\|_2^2\le\sigma^2$ for all
$1\le u\le r_1+r_2$. The union of these sets and functions yield a
super-regular $(r_1+r_2,\sigma^2)$ dominating system of the function
$f\circ g$ consisting of the sets $B_1,\dots,B_{r_1+r_2}$ and the
functions $h_1,\dots,h_{r_1+r_2}$.
 
We shall prove Lemma 4 with the help of the above observation,
formula (4.25) and some properties of the operators $P_j$ and
$R_{j,j'}$. First we show by induction that the function
$\prodd_{t=1}^p P_{j_t} R_{j_t,j'_t} f\circ g(x_1,\dots,x_{k_1+k_2})
=\overline{f\circ g}_{B(\Cal N_1,\Cal N_1)}$ also has a
super-regular $(r_1+r_2,\sigma^2)$ dominating system.
 
We shall prove the following statement which implies the above
relation. Let a function $G(x_{s_1},\dots,x_{s_m})$ have an
$(r,\sigma^2)$ super-regular dominating system and take two
indices $1\le j,\le k_1<j'\le k_1+k_2$  such that
$\{j,j'\}\subset\{s_1,\dots,s_m\}$. Then the function $P_jR_{j,j'}G$
also has an $(r,\sigma^2)$ super-regular dominating system. This
statement enables us to show with the help of a simple induction that
together with the function $f\circ g$ also the function
$\overline{f\circ g}_{B(\Cal N_1,\Cal N_1)}=\prodd_{t=1}^p
P_{j_t}R_{j_t,j'_t}f\circ g$ has a super-regular
$(r_1+r_2,\sigma^2)$ dominating system.
 
To prove the above statement let us consider a super-regular
$(r,\sigma^2)$ dominating system of the function $G$ consisting of
a partition $B_1,\dots,B_r$ of the set $\{s_1,\dots,s_m\}$ and some
functions $h_l$, $1\le l\le r$, depending on the arguments indexed
by the set $B_l$. To construct a super-regular dominating system of
the function $P_jR_{j,j'}G$ let us introduce the following notations.
Let $u_1$ and $u_2$ denote those indices for which $j\in B_{u_1}$ and
$j'\in B_{u_2}$, and let us denote the elements of these sets as
$B_{u_1}=\{j,v(1),\dots,v(m)\}$, $B_{u_2}=\{j',v'(1),\dots,v'(m')\}$.
Then we define $\bar B_u=B_u$, $\bar h_u=h_u$ if $u\neq u_1$ and
$u\neq u_2$, $\bar B_{u_1}=B_{u_1}\setminus \{j\}$,
$\bar B_{u_2}=B_{u_2}\setminus \{j'\}$, and
$$
\align
\bar h_{u_1}(x_{v(1)},\dots,x_{v(m)})&=
\(\int h^2_{u_1}(z_j,x_{v(1)},\dots,x_{v(m)})\mu(\,dz_j)\)^{1/2},\\
\bar h_{u_2}(x_{v'(1)},\dots,x_{v'(m')})&=
\(\int h^2_{u_2}(z_{j'},x_{v'(1)},\dots,x_{v'(m')})
\mu(\,dz_{j'})\)^{1/2}
\endalign
$$
if $u=u_1$ or $u=u_2$. We claim that the above defined sets $\bar
B_l$ and functions $\bar h_l$, $1\le l\le r$, supply a
super-regular $(r,\sigma^2)$ dominating system of the function
$P_jR_{j,j'}G$. To show this let us first observe that because
of the Schwarz inequality
$$
\aligned
&|P_jR_{j,j'} h_{u_1}(x_j,x_{v(1)},\dots,x_{v(m)})
h_{u_2}(x_{j'},x_{v'(1)},\dots,x_{v'(m')})|\\
&\qquad \le \bar h_{u_1}(x_{v(1)},\dots,x_{v(m)})
\bar h_{u_2}(x_{v'(1)},\dots,x_{v'(m')}).
\endaligned\tag4.27
$$
Since the product of the functions $h_l$, $1\le l\le r$, gives an
upper bound of the function $G$, formula (4.27) and the definition
of the functions $\bar h_l$ imply that the product of the functions
$\bar h_l$ yields an upper bound for the function $P_jR_{j,j'}G$.
The functions $\bar h_{u_j}$, $j=1,2$, satisfy also the relation
$\|\bar h_{u_j}\|^2=\|h_{u_j}\|^2\le\sigma^2$. It is not
difficult to see with the help of these relations that the sets
$\bar B_l$ and $\bar h_l$, $1\le l\le r$, yield a super-regular
$(r,\sigma^2)$ dominating system of the function $R_{j.j'}P_j G$.
 
Thus we have proved that, with the notation  (4.24) and (4.26), the
function $\overline{f\circ g}_{B(\Cal N_1,\Cal N_1)}
(x_{s(1)},\dots,x_{s(k_1+k_2-2p)})$ together with the diagram $B(\Cal
M)$ attached to it, where $\Cal M=\Cal N\setminus \Cal N_1=
\{(j_{v(1)},j_{v(1)}'),\dots,(j_{v(l-p)},j_{v(l-p)}')\}$, and the
indices of the arguments of the function
$\overline{f\circ g}_{B(\Cal N_1,\Cal N_1)}$ are
$$
\align
\{s(1),\dots,s(k_1-p)\}&=\{1,\dots,k_1\}\setminus
\{j_{u(1)},\dots,j_{u(p)}\},\\
\{s(k_1-p+1),\dots,s(k_1+k_2-2p)\}&=\{k_1+1,\dots,k_1+k_2\}
\setminus \{j_{u(1)}',\dots,j_{u(p)}'\},
\endalign
$$
has a regular $(r_1+r_2,\sigma^2)$ dominating system. We want to
prove that the function
$$
\overline{f\circ g}_{B(\Cal N,\Cal N_1)}
=\prodd_{m=1}^{l-p} R_{j_{v(m)},j_{v(m)}'} \overline{f\circ
g}_{B(\Cal N_1,\Cal N_1)}        \tag4.28
$$
has an $(r_1+r_2-(l-p),\sigma^2)$ dominating system. Let us recall
that $|\Cal M|=l-p$ and the operators $R_{j_{v(m)},j_{v(m)}'}$
considered in formula (4.28) are commutative.
 
The above statement will be proved by the successive application of
a statement formulated below. Let a function $G(x_{s(1)},\dots,
x_{s(m)})$ and a diagram $B(\Cal M)$ attached to it
have a regular $(r,\sigma^2)$ dominating system consisting of some
partition $B_1,\dots,B_r$ of the set $\{s(1),\dots,s(m)\}$ and
some functions $h_l$, $1\le l\le r$, depending on variables indexed
by the set $B_l$. Let us consider two sets $B_l$ and $B_{l'}$ such
that there are edges of the diagram $\Cal M$ whose end-points are in
the sets $B_l$ and $B_{l'}$. By a new enumeration of the sets of the
partition if it is needed we may assume that $l=r-1$ and $l'=r$. Let
$\Cal M_0=\{(s(1),s(1)'),\dots,(s(q),s(q)')\}$, $q\ge 1$, denote
the set of all edges of $\Cal M$ which connect some points of the
sets $B_{r-1}$ and $B_r$. Here we do not specify which element of
a pair $(s(t),s(t)')$ is in the set $B_{r-1}$ and which element
is in the set $B_r$, $1\le t\le q$. We claim that the function
$\prodd_{t=1}^q R_{s(t),s(t)'}G$ together with the diagram $B(\Cal
M\setminus\Cal M_0)$ has a regular $(r-1,\sigma^2)$ dominating
system. (Let us recall that the diagram attached to the function
$G$ has the role to list the pairs $(j,j')$ for which we want to
apply the operator $R_{j,j'}$.) More explicitly, we can get a
regular $(r-1,\sigma^2)$ dominating system $(\bar B_l,\bar h_l)$,
$1\le l\le r-1$, in the following way. Put $\bar B_l=B_l$ and $\bar
h_l=h_l$, $1\le l\le r-2$, $\bar B_{r-1}=B_{r-1}\cup B_r\setminus
\{s(1)',\dots,s(q)'\}$ and $\bar h_{r-1}=\prodd_{t=1}^q
R_{s(t),s(t)'}h_{r-1}h_r$. It is rather straightforward to check
that the above sets $\bar B_l$ and functions $\bar h_l$, $1\le l\le
r-1$ constitute a regular $(r-1,\sigma^2)$ dominating system. The
only statement which demands some special consideration is the
inequality $\|\bar h_{r-1}\|_2^2\le\sigma^2$. We can see with
the help of the Schwarz inequality (similarly to the
corresponding argument applied in the proof of formula (3.14))
that
$$
\|\bar h_{r-1}\|_2^2=\left\|\prodd_{t=1}^q
R_{s(t),s(t)'}h_{r-1}h_r\right\|_2^2\le
\|h_{r-1}\|_2\|h_r\|_2\le\sigma^2.
$$
 
By applying the above fact for the function at the right-hand side
of (4.28) $s\le l-p$ times we get that the function
$\overline{f\circ g}_{B(\Cal N,\Cal N_1)}$ is $(r_1+r_2-s,\sigma^2)$
dominated with some $s\le l-p$. (At each step we take into
consideration at least one new operator $R_{j_{v(m)},j_{v(m)}'}$ in
formula (4.28), so the procedure finishes in $s\le l-p$ steps.)
This property implies that the function
$\overline{f\circ g}_{B(\Cal N,\Cal N_1)}$ has an $(r_1+r_2-s,
\sigma^2)$ dominating system with some number $s\le l-p$. To give a
complete proof for the main result of Lemma~4 let us observe that
in the case of a `too good estimate' when $s<l-p$ we can construct
an $(r,\sigma^2)$ dominating system of the function
$\overline{f\circ g}_{B(\Cal N,\Cal N_1)}$ with
$r=r_1+r_2-(l-p)$. We get such a system by taking the union of
some sets $B_l$ of an $(\bar r,\sigma^2)$ dominating system if $\bar
r>r$, and by considering the product of the functions $h_l$
corresponding to them.
 
Finally, the last statement of Lemma 4 is a simple consequence of
formula (3.14). Indeed, if $f$ is $(r_1,\sigma^2)$ and $g$ is
$(r_2,\sigma^2)$ dominated, then $\|f\|_2^2\le \sigma^{2r_1}$,
$\|g\|_2^2\le \sigma^{2r_2}$, hence $\|f\circ g_{B(\Cal N,\Cal
N_1)}\|_2^2\le\sigma^{r_1+r_2}$ by relation (3.14). On the other
hand, $\|f\|_\infty\le1$ and $\|g\|_\infty\le1$, hence
$\|f\circ g_{B(\Cal N,\Cal N_1)}\|_\infty\le1$. These inequalities
mean that $f\circ g_{B(\Cal N,\Cal N_1)}$ is an
$(1,\sigma^{(r_1+r_2)/2})$ dominated function as we claimed.
 
\beginsection 5. The proof of Proposition 2 together with a
technical lemma
 
First I prove Proposition 2 with the help of Lemma~5 formulated in
Section~2. Then I also prove Lemma~5.
\medskip\noindent
{\it The proof of Proposition 2 with the help of Lemma~5.}\/ By
Lemma~3 formula (3.17) holds for $m=0$ and all $k=1,2,\dots$ with
the choice $C(k,0)=C^k$ with some $C>0$. Let us assume that it holds
for some $m$ and all $k\ge1$ with appropriate norming constants
$C(k,m)$. Then I prove it by induction for $m+1$ with
appropriate norming constants $C(k,m+1)$.
 
To do this let us fix some function $f(x_1,\dots,x_k)$ of $k$
variables, consider the functions $f\circ f_{l,p}$ defined in
formula (3.12) and give a good bound for the $2M=2^{m+1}$-th
moment of the random integrals $J_{n,2k-l-p}(f\circ f_{l,p})$.
Then the $4M=2^{m+2}$-th moment of the random integral $J_{n,k}(f)$
(in the case $n\ge 2kM$) can be estimated by means of these bounds
and the diagram formula.
 
Let $f=f(x_1,\dots,x_k)$ be a function of $k$ variables, $k\ge1$,
which is $\(r,\sigma^{2/r}\)$ dominated with some $r\ge1$. I show
with the help of Lemma 4 and formula (3.17) with the parameter
$M=2^m$ that the inequality
$$
\aligned
&E \(n^{-(l-p)/2}J _{n,2k-l-p} (f\circ f_{l,p})\)^{2M} \le
\(\frac{C(2k-l-p,m)}{(2k-l-p)^{2k-l-p}}\)^M\((2M)^k
\sigma^2\)^{2M} \\
&\qquad\cdot\frac{2^{-2kM}}{(2k)^{(l-p)M} M^{2lM} }
\max \(1,\(\frac{2kM}{n\sigma^{2/r}}
\)^{2M\min(k,r)}\), \quad 0\le p\le l\le k
\endaligned\tag5.1
$$
holds if $2kM\le n$. This formula also holds in the degenerated case
$p=l=k$ with the choice $C(0,m)=1$. Indeed, in this case $2k-p-l=0$,
hence $f\circ f_{k,k}$ is a constant, and $E J_{n,0}(f\circ
f_{k,k})^{2M}=(f\circ f_{k,k})^{2M}\le \(\sigma^2\)^{2M}$. The last
inequality follows from the definition of the function $f\circ
f_{k,k}=f\circ f_{B(\Cal N,\Cal N)}$, $|\Cal N|=k$, given in formula
(3.12) and the estimate $|f\circ f_{B(\Cal N,\Cal N_1)}|
\le\|f\|_2^2\le\sigma^2$ if $\Cal N=\Cal N_1$, $|\Cal N|=k$, ($k$ is
the number of variables of the function $f$). The last inequality
follows from relation (3.15). These relations imply that $|f\circ
f_{k,k}|\le \sigma^2$. We had to consider this case separately,
because our inductive hypothesis (3.17) is formulated only for
functions with $k\ge1$ variables.
 
The integral $J_{n,2k-l-p}(f\circ f_{l,p})$ is the average of the
random integrals $J_{n,2k-l-p}(f\circ f_{B(\Cal N,\Cal N_1)})$ with
the functions $f_{B(\Cal N,\Cal N_1)}(x_1,\dots,x_{2k-l-p})$ such
that $B(\Cal N,\Cal N_1)\in\Cal B(l,p)$. Hence the triangular
inequality for the $L_{2M}$ norm in the probability space $(\Omega,
\Cal A,P)$ where the random variables $\xi_j$ are living implies
that it is enough to prove such a version of formula (5.1) where the
random integral $J_{n,k}(f\circ f_{l,p})$ at the left-hand side is
replaced by an arbitrary random integral $J_{n,2k-l-p}(f\circ
f_{B(\Cal N,\Cal N_1)})$ such that $B(\Cal N,\Cal N_1)\in \Cal
B(l,p)$.
 
Let us apply formula (3.17) for a random integral
$J_{n,2k-l-p}(f\circ f_{B(\Cal N,\Cal N_1)})$ such that $B(\Cal N,
\Cal N_1)\in \Cal B(l,p)$. (In this case $(2k-l-p)M\le n$, and our
inductive hypothesis allows the application of formula (3.17).)
Let us consider the cases $2r\ge l-p+1$ and $2r\le l-p$ separately.
In the first case the function $f\circ f_{B(\Cal N, \Cal N_1)}$ is
$\(2r-l+p,\sigma^{2/r}\)$ or in other notation
$\(2r-l+p,\bar\sigma^{2/(2r-l+p)}\)$ dominated with
$\bar\sigma=\sigma^{(2r-l+p)/r}$ by Lemma~4. Let us apply
relation (3.17) with the choice $\bar k=2k-l-p$, $\bar r=2r-l+p$ and
$\bar\sigma=\sigma^{(2r-l+p)/r}$ in this case. Observe that
$\sigma^{2/r}=\bar\sigma^{2/\bar r}$. We get that
$$
\align
&E \(n^{-(l-p)/2}J_{n,2k-l-p}(f\circ f_{B(\Cal N,\Cal N_1)})\)^{2M}
\le \(\frac{C(2k-l-p,m)}{(2k-l-p)^{2k-l-p}}\)^M\((2M)^k
\sigma^2\)^{2M}\\
&\qquad\qquad
\cdot\frac{2^{-2kM}}{M^{(l+p)M}}\frac1{(n\sigma^{2/r})^{(l-p)M}}
\max \(1,\(\frac{2kM}{n\sigma^{2/r}}
\)^{M\min(2k-l-p,2r-l+p)}\) \\
&\qquad=\(\frac{C(2k-l-p,m)}{(2k-l-p)^{2k-l-p}}\)^M\((2M)^k
\sigma^2\)^{2M}\\
&\qquad\qquad\cdot\frac{2^{-2kM}}{M^{2lM}(2k)^{(l-p)M}}
\max\(\(\frac{2kM}{(n\sigma^{2/r}}\)^{(l-p)M},
\(\frac{2kM}{n\sigma^{2/r}}\)^{M\min(2k-2p,2r)}\) \\
&\qquad\le\(\frac{C(2k-l-p,m)}{(2k-l-p)^{2k-l-p}}\)^M\((2M)^k
\sigma^2\)^{2M}\\
&\qquad\qquad\cdot\frac{2^{-2kM}}{M^{2lM}(2k)^{(l-p)M}}
\max\(1,\(\frac{2kM}{n\sigma^{2/r}}\)^{2M\min(k,r)}\),
\tag5.2
\endalign
$$
since the exponents $(l-p)M$ and $\min(2k-2p, 2r)M$ in the last but
one inequality of formula (5.2) satisfy the inequality $0\le (l-p)M\le
\min 2M(k,r)$, and $0\le\min(2k-2p, 2r)M\le \min 2M(k,r)$.
 
In the second case when $2r\le l-p$ we know by Lemma~4 that the
function $f\circ f_{B(\Cal N,\Cal N_1)}$ is $(1,\sigma^2)$ dominated.
In this case we apply formula (3.17) with the $\bar k=2k-l-p$, $\bar
r=1$ and $\bar\sigma=\sigma$. We get that
$$ \allowdisplaybreaks
\align
&E\(n^{-(l-p)/2}J_{n,2k-l-p}(f_{B(\Cal N,\Cal N_1)})\)^{2M} \le
\(\frac{C(2k-l-p,m)}{(2k-l-p)^{2k-l-p}}\)^M\((2M)^k
\sigma^2\)^{2M}\\
&\qquad\qquad
\cdot\frac{2^{-2kM}}{M^{(l+p)M}}\frac1{\sigma^{2M}n^{(l-p)M}}
\max \(1,\(\frac{2kM}{n\sigma^2} \)^M\) \\
&\qquad=\(\frac{C(2k-l-p,m)}{(2k-l-p)^{2k-l-p}}\)^M\((2M)^k
\sigma^2\)^{2M}  \tag5.3   \\
&\qquad\qquad
\cdot\frac{2^{-2kM}}{M^{(l+p)M}}\max\(\(\frac1{n^{(l-p)}
\sigma^{2}}\)^M ,\(\frac{2kM}{n^{l-p+1}\sigma^4}\)^M\)
\endalign
$$
I claim that
$$
\max\(\frac1{n^{(l-p)}\sigma^2},
\frac{2kM}{n^{l-p+1}\sigma^4}\)\le\dfrac1{(2kM)^{l-p}} \max\(1,\(\frac
{2kM}{n\sigma^{2/r}}\)^{2r}\). \tag5.4
$$
Relation (5.4) together with formula (5.3) and the relation $k\ge r$
if $2r\le l-p$ yield that in this case also the inequality
$$
\align
&E \(n^{-(l-p)/2} J_{n,2k-l-p}(f\circ f_{B(\Cal N,\Cal N_1)})\)^{2M}
\le\(\frac{C(2k-l-p,m)}{(2k-l-p)^{2k-l-p}}\)^M\((2M)^k
\sigma^2\)^{2M}\\
&\qquad\qquad \cdot\frac{2^{-2kM}}{M^{2lM}}
\dfrac1{(2k)^{(l-p)M}} \max\(1,\(\frac
{2kM}{n\sigma^{2/r}}\)^{2M\min(r,k)}\) \tag5.5
\endalign
$$
holds. Relations (5.2) and (5.5) give a bound on
$E\(n^{-(l-p)/2}J_{n,2k-l-p}(f_{B(\Cal N,\Cal N_1)})\)^{2M}$
in both cases $2r\ge l-p+1$ and $2r\le l-p$, thus they imply formula
(5.1). Hence to complete the proof of formula (5.1) it is enough to
check relation (5.4).
 
Actually I prove a stronger inequality than formula (5.4). I show
that if the inequality $2kM\le n$ holds, and this was assumed among
the conditions of relation (5.1), then the left-hand side of (5.4)
can be bounded by
$\dfrac1{(2kM)^{l-p}}\(\dfrac{2kM}{n\sigma^{2/r}}\)^{2r}$, i.e.\ by
the second term of the maximum at the right-hand side of (5.4). This
follows from the following estimations:
$$
\frac1{n^{l-p}\sigma^2}
=\dfrac1{(2kM)^{l-p}}\(\frac{2kM}{n\sigma^{2/r}}\)^{2r}
\(\frac n{2kM}\)^{2r-(l-p)}\sigma^2\le
\dfrac1{(2kM)^{l-p}}\(\frac{2kM}{n\sigma^{2/r}}\)^{2r}
$$
and
$$
\frac{2kM}{n^{l-p+1}\sigma^4}=\(\frac{2kM}{n\sigma^{2/r}}\)^{2r}
\frac1{(2kM)^{l-p}} \(\frac n{2kM}\)^{2r-(l-p)-1}
\le\(\frac{2kM}{n\sigma^{2/r}}\)^{2r}
\frac1{(2kM)^{l-p}},
$$
since $2kM\le n$, $2r-(l-p)\le0$ and $\sigma^2\le1$.
 
The diagram formula (3.11) together with the triangular inequality
in the $L_{2M}$ norm yield that
$$
\aligned
&E(J_{n,k}(f)^{4M}=
E\(\sum_{l=0}^k\sum_{p=0}^l
\frac{(2k-l-p)!}{(k-l)!^2(l-p)!p!} n^{-(l-p)/2}
\cdot J_{n,2k-l-p}(f\circ f_{l,p})\)^{2M}\\
&\qquad \le \(\sum_{l=0}^k\sum_{p=0}^l
\frac{(2k-l-p)!}{(k-l)!^2(l-p)!p!} \(E\(n^{-(l-p)/2}
J_{n,2k-l-p}(f\circ f_{l,p}\)^{2M} \)^{1/2M}\)^{2M}.
\endaligned\tag5.6
$$
Relations (5.1) and (5.6) imply the inequality
$$
\aligned
&E(J_{n,k}(f)^{4M}\le \((2M)^k\sigma^2\)^{2M}
\max \(1,\(\frac{2kM}{n\sigma^{2/r}}\)^{2M\min(k,r)}\)\\
&\qquad\cdot\(\sum_{l=0}^k\sum_{p=0}^l
\dfrac{(2k-l-p)!}{(k-l)!^2(l-p)!p!}
\(\frac{C(2k-l-p,m)}{(2k-l-p)^{2k-l-p}}\)^{1/2}
\frac{(2k)^{-(l-p)/2}} {2^k\cdot M^l}\)^{2M}
\endaligned\tag5.7
$$
for the random integral $J_{n,k}(f)$ of an $\(r,\sigma^{2/r}\)$
dominated function $f$.
By writing the term $M^l$ in the form $2^{ml}$ in (5.7) and
comparing formulas (3.17) and (5.7) we can see that relation
(3.17) also holds for the new parameter $m+1$ with some new
(positive) constants $C(k,m+1)$, $k=1,2,\dots$, if we can choose
them in such a way that they satisfy the inequalities
$$
\dfrac{C(k,m+1)}{k^k}\ge
\sum_{l=0}^k\sum_{p=0}^l
\dfrac{(2k-l-p)!}{(k-l)!^2(l-p)!p!}
\(\frac{C(2k-l-p,m)}{(2k-l-p)^{2k-l-p}}\)^{1/2}
\frac{(2k)^{-(l-p)/2}}{2^k\cdot 2^{lm}} \tag5.8
$$
for all $k=1,2,\dots$. We show with the help of Lemma~5 that it
is possible to choose such constants $C(k,m)$, $k=0,1,\dots$,
$m=0,1,\dots$ which satisfy both relation (5.8) and the
inequalities $\supp_{0\le m<\infty} C(k,m)\le C_k<\infty$ for all
$k=0,1,2,\dots$, and beside this also the inequalities $C(0,m)\ge1$,
$m=1,2,\dots$, and $C(k,0)\ge C^k$, $k=0,1,2,\dots$, hold for all
$k=0,1,\dots$ with a prescribed number $C>0$. (The last two
inequalities are needed to satisfy relation (5.1) in the degenerate
case $k=0$ and relation (3.17) for the starting step $m=0$.) The
proof of the above statement completes the proof of Proposition~2.
 
To find a sequence $C(n,k)$ satisfying relation (5.8) first I show
that a sequence $\bar C(n,k)$ satisfying Lemma~5 also satisfies
the system of inequalities
$$
\aligned
\bar C(k,m+1)^2&\ge\(\dfrac{(2k-l-p)!}{(k-l)!^2(l-p)!p!}\)^2
\frac{k^{2k}}{(2k-l-p)^{2k-l-p}} \frac{(2k)^{-(l-p)}}{2^{2k}\cdot
2^{2lm}} \bar C(2k-l-p,m) \\
&\qquad \text{for all } 0\le k<\infty,
\text { and } \; 0\le m<\infty \text { and } 0\le p<l\le k.
\endaligned \tag5.9
$$
(Formula (5.9) is a weakened version of (5.8)  where it is demanded
that if the double sum at the right-hand side of (5.8) is replaced
by the maximum of these terms, then this new expression should be
smaller than the left-hand side.)
 
Indeed, to show that a sequence satisfying Lemma~5 also satisfies
(5.9) observe that by the estimate $\binom nm\le 2^n$ for all $0\le
m\le n$ and the rather rough bound  $2^{2l} l!\ge l^l$
$$
\align
\dfrac{(2k-l-p)!}{(k-l)!^2(l-p)!p!}&\le (2k-l-p)^{l-p}\binom
{2k-2l}{k-l}\frac1{l!}\binom lp\\
&\le (2k-l-p)^{l-p} 2^{2(k-l)}\frac{2^{2l}}{l^l} 2^l
=2^{2l} \frac{2^{2k}(2k-l-p)^{l-p}}{(2l)^l},
\endalign
$$
and this estimation implies that inequality (5.9) follows from the
relation
$$
\bar C(k,m+1)^2\ge
\frac{1}{2^{2lm}} \frac{2^{4l}(2k)^{2k-l+p}(2k-l-p)^{3l-p-2k}}
{(2l)^{2l}} \bar C(2k-l-p,m)
$$
i.e.\ from inequality (3.18) imposed in Lemma~5.
 
We claim that the numbers $C(k,m)=AB^k \bar C(k,m)$ defined  with
the help of the above numbers $\bar C(k,m)$ and sufficiently large
constants $A>1$ and $B>1$ satisfy formula (5.8) together with the
additional properties we have imposed. Indeed, we get with such
a choice that
$$
\align
&\sum_{l=0}^k\sum_{p=0}^l
\dfrac{(2k-l-p)!}{(k-l)!^2(l-p)!p!}
\(\frac{C(2k-l-p,m)}{(2k-l-p)^{2k-l-p}}\)^{1/2}
\frac{(2k)^{-(l-p)/2}}{2^k\cdot 2^{lm}}\\
&\qquad\le\frac1{k^k} \sum_{l=0}^k\sum_{p=0}^l
\[A B^{2k-l-p}\bar C(k,m+1)^2\]^{1/2}
= \frac{C(k,m+1)}{\sqrt A k^k}\sum_{l=0}^k\sum_{p=0}^l B^{-(l+p)/2}\\
&\qquad\le\frac {C(k,m+1)}{\sqrt A k^k}\(\sum_{j=1}^\infty B^{-j/2}\)^2
\le\frac{C(k,m+1)}{k^k}
\endalign
$$
if we choose the number $A=A(B)$ sufficiently large. Because of
Lemma~5 the numbers $C(k,m)$ also satisfy the relation $\supp_{0\le
m<\infty}C(k,m)\le C_k<\infty$ with an appropriate constant $C_k$.
Beside this, if we choose $B\ge C$, then because of the relation
$\bar C(k,0)=1$ formulated in Lemma~5 also the relation $C(k,0)\ge
C^k$ holds for all $k=0,1,\dots$. The relation $C(0,m)=A\bar
C(0,m)=A\ge1$ also holds. The proof of Proposition 2 with the help
of Lemma 5 is finished. \medskip\noindent
{\it The proof of Lemma 5.}\/ Let us introduce the numbers
$$
A(l,p,k,m)=2^{2l(4-m)} \frac{(2k)^{2k-l+p}(2k-l-p)^{3l-p-2k}}
{(2l)^{2l}}, \quad 0\le p\le l\le k,\; m\ge0.
$$
To construct a sequence $\bar C(k,m)$ which satisfies Lemma~5 it
is enough to show that there exist such numbers $D(m)\ge 1$,
$m=0,1,\dots$, which satisfy the inequalities $D(m)^{2k}\ge
A(l,p,k,m)$ for all $0\le p\le l\le k$, $m=0,1,\dots$, and
$\prodd_{m=0}^\infty D(m)<\infty$. Indeed, if there exists such a
sequence $D(m)$, then the numbers $\bar C(k,0)=1$, $\bar
C(k,m)=\(\prodd_{p=0}^{m-1} D(p)\)^k$ for $m\ge1$, $k\ge0$, satisfy
Lemma~5, since the above defined numbers satisfy the inequality
$$
\bar C(k,m+1)^2\ge \bar C(2k-l-p,m) D(m)^{2k}\ge
A(l,p,k,m) \bar C(2k-l-p,m)
$$
for all $0\le p\le l\le k$.
 
To show that we can define numbers $D(m)$ with the above required
properties let us consider the extension of the sequence
$A(l,p,k,m)$ to the real valued function
$$
A_{k,m}(u,v)=2^{2(4-m)v} \frac{(2k)^{2k+u-v} {(2k-u-v)}^{3v-u-2k}}
{(2v)^{2v}},\quad 0\le u\le v\le k,
$$
and let us estimate its maximum. To get a simpler estimate let us
bound the function $A_{k,m}(u,v)$ in the following way. Observe that
$$
\align
A_{k,m}(u,v)&\le \(\frac{2k}{2k-u-v}\)^{v-u} A_{k,m}(u,v)=2^{2(4-m)v}
\frac{(2k)^{2k}}{(2k-u-v)^{2k-2v}(2v)^{2v}} \\
&\le 2^{2(4-m)v}\frac{(2k)^{2k}}{(2k-2v)^{2k-2v}(2v)^{2v}}
\endalign
$$
for $0\le u\le v\le k$, and this inequality implies that
$$
\sup_{0\le u\le v\le k} A_{k,m}(u,v)\le \sup_{0\le v\le k}B_{k,m}(v)
$$
with
$$
B_{k,m}(v)=2^{2(4-m)v} \frac{(2k)^{2k}}{(2k-2v)^{2k-2v}(2v)^{2v}}.
$$
Simple differentiation of the function $\log B_{k,m}(v)$ shows that
the function $B_{k,m}(v)$ takes its maximum in the point
$\bar v=\dfrac k{2^{(m-4)}+1}$, and $B_{k,m}(\bar v)
=\(1+2^{4-m}\)^{2k}$.
 
The above calculations imply that the numbers $D(m)=1+2^{4-m}$,
$m=0,1,2,\dots$, satisfy the inequality $\supp_{0\le p\le l\le
k}A(l,p,k,m)\le D(m)^{2k}$. The relation $\prodd_{m=0}^\infty
D(m)<\infty$ also holds. Hence, as the argument at the beginning
of the proof shows the numbers  $\bar C(k,0)=1$, $\bar
C(k,m)=\(\prodd_{p=0}^{m-1} D(p)\)^k$ for $m\ge1$ satisfy the
properties demanded in Lemma~5.
 
\bigskip\noindent
{\bf References:} \medskip
 
\item{1.)} Arcones, M. A. and Gin\'e, E. (1993) Limit theorems for
$U$-processes. {\it Ann. Probab.} {\bf 21}, 1494--1542
\item{2.)} Arcones, M. A. and Gin\'e, E. (1994) $U$-processes
indexed by Vapnik--\v{C}ervonenkis classes of functions with
application to asymptotics and bootstrap of $U$-statistics with
estimated parameters. {\it Stoch. Proc. Appl.}  {\bf 52}, 17--38
\item{3.)} Major, P. (1981) Multiple Wiener--It\^o integrals. {\it
Lecture Notes in Mathematics\/} {\bf 849}, Springer--Verlag, Berlin
Heidelberg, New York,
\item{4.)} Major, P. (1988) On the tail behaviour of the distribution
function of multiple stochastic integrals. {\it Probability Theory
and Related Fields}, {\bf 78},  419--435
\item{5.)} Major, P. and Rejt\H{o}, L. (1988) Strong embedding of
the distribution function under random censorship, {\it Annals of
Statistics}, {\bf 16}, 1113--1132
\item{6.)} Major, P. and Rejt\H{o}, L. (1998) A note on nonparametric
estimations. In the conference volume to the 65. birthday of Mikl\'os
Cs\"org\H{o}.759--774
\item{7.)} Major, P. (2003) An estimate on the maximum of a nice class
of stochastic integrals. submitted to {\it Probability Theory
and Related Fields},
\item{8.)} Surgailis, D. (1984) On multiple Poisson stochastic
integrals and associated Markov semigroups. {\it Probab. Math.
Statist.} 3. no. {\bf 2} 217-239
 
\bigskip\bigskip\bigskip
Supported by the OTKA foundation Nr. 037886

\bye